\documentclass[12pt,leqno]{amsart}

\usepackage{amsmath}
\usepackage{amscd}
\usepackage{amssymb}
\usepackage{latexsym}
\makeatletter

\@addtoreset{equation}{thm}
\makeatother

\newtheorem{thm}{Theorem}[section]
\newtheorem{pr}[thm]{Proposition}
\newtheorem{df}[thm]{Definition}
\newtheorem{lm}[thm]{Lemma}
\newtheorem{cor}[thm]{Corollary}
\newtheorem{cn}[thm]{Conjecture}
\newtheorem{rmk}[thm]{Remark}
\newtheorem{ex}[thm]{Example}

\def\disc{\text{\rm disc }}

\def\Spec{\text{\rm Spec }}
\def\Gal{\text{\rm Gal}}

\input xy \xyoption{all} \CompileMatrices

\begin{document}

\title{The second Stiefel-Whitney class
of $\ell$-adic cohomology}
\author{{\sc Takeshi Saito}}
\address{Department of Mathematical Sciences, 
University of Tokyo, 
Tokyo 153-8914, Japan}
\email{t-saito@ms.u-tokyo.ac.jp}
\maketitle

\begin{abstract}
For a proper smooth variety
of even dimension
over a field
of characteristic different from 2 or $\ell$,
the second Stiefel-Whitney class
of the $\ell$-adic cohomology
and the second Hasse-Witt
class of the de Rham cohomology
are both defined in the second
Galois cohomology.
We state a conjecture on
their relation
and give several evidences.
\end{abstract}

The cohomology
of middle degree
of a proper smooth
variety
of even dimension
carries a non-degenerate
symmetric bilinear form.
This defines
various Stiefel-Whitney
classes
depending on the cohomology
theory we consider.
The second
Stiefel-Whitney
class of $\ell$-adic
cohomology
allows us to deduce
that the $L$-function
has the positive sign
in the functional equation,
from a reasonable
hypothesis \cite{S2}.
For the de Rham cohomology,
we also call the
Stiefel-Whitney class
the Hasse-Witt class
since it is defined
by a quadratic form.

In this paper,
we propose a conjecture
Conjecture \ref{cn}
comparing
the Stiefel-Whitney
classes of $\ell$-adic
cohomology
and of de Rham cohomology
and prove several 
evidences.
In the case where the variety
is of dimension $0$,
it amounts to
a formula of Serre
on the Hasse-Witt
invariant of the trace form
for a finite separable
extension
\cite{Se}.
Conjecture \ref{cn}
may also be 
regarded as 
a version in
degree 2
of the formula
for the determinant
of cohomology proved
in \cite{S1}.

We state the main
Conjecture \ref{cn}
in Section \ref{scn}
after some preliminaries
on the Stiefel-Whitney
classes of orthogonal
Galois representations
and of quadratic
vector spaces in
Section \ref{scls}.
We give an apparently simpler
reformulation in Corollary \ref{corg}
by using a graded variant.

As evidences
for Conjecture \ref{cn},
we prove the following
cases in
Theorem \ref{thm1}
under some
additional conditions:
1. The base field 
is a finite extension
of ${\mathbb Q}_p$
for $p\neq \ell$
and the variety
has a certain mild degeneration.
2. The base field
is a finite unramified
extension of 
${\mathbb Q}_\ell$
and the variety
has a good reduction.
3. The base field is 
${\mathbb R}$.
4. The base field
is an extension of
$\bar {\mathbb Q}$.
5. The base field
is arbitrary
and the variety
is a smooth
hypersurface in 
a projective space.

Curiously,
Theorem \ref{thm1}
implies that
the second Stiefel-Whitney class
of $\ell$-adic cohomology
in fact may depend on $\ell$,
as in Example \ref{ex},
contrary to
the first Stiefel-Whitney class
which is independent of $\ell$
as a consequence of the Weil conjecture
(see Lemma \ref{lm1}).

We formulate a 
generalization
Conjecture \ref{cn2}
of Conjecture \ref{cn} 
for families
in Section \ref{sfam}
after introducing
the Stiefel-Whitney classes
for symmetric complexes
in Section \ref{sswcx}.
A similar construction
is studied by Balmer
in \cite{Ba}.
Contrary to there,
we are more interested in the invariants
of individual complexes
rather than the invariants of
the categories.
The author learned from \cite{CNET}
that a similar construction
is also studied in \cite{Wa}.
The rest of the article
is devoted to
the proof of Theorem \ref{thm1}.

In Section \ref{sdeg},
we prove the assertion 1
of Theorem \ref{thm1}
using nearby cycles
and the sheaves
of differential forms
with logarithmic poles.
We observe in Lemma \ref{lmdbl}
that the mysterious appearance
of the term involving $2$
in Serre's formula
arises from the multiplicity
of the exceptional divisor
of the blow-up
at an ordinary double point.
We verify
that the assertion 2 is
essentially proved
in \cite{S2}
in Section \ref{scris}.
In Section \ref{sHdg},
we prove the assertion 3
using Hodge structures.
We prove the assertion 4
by partly proving the
generalization
Conjecture \ref{cn2}
for families by
a transcendental argument
in Section \ref{strs}.
In Section \ref{sart},
using the universal family
of hypersurfaces,
we prove the assertion 5
by a global arithmetic argument.
In the proof of 5.,
we combine all the results
obtained in
the other parts of the
article.

The author thanks Paul Balmer
for discussion
on symmetric complexes.
The author also thanks
Asher Auel for his interest
which encouraged him
to complete the article.
The author thanks Luc Illusie
for informing an unpublished preprint \cite{Su}.

The research is partly
supported by
Grants-in-aid for
Scientific Research 
S-19104001.

\section{Stiefel-Whitney classes
and Hasse-Witt classes}\label{scls}

We recall some basic definitions
to formulate a conjecture.
In order to distinguish
the Stiefel-Whitney
classes of
orthogonal representations
from those
of quadratic forms,
we will write
$sw$ for the
former and
$hw$ for the latter.
We call the latter
the Hasse-Witt class.

Let $\pi$ be
a group
and $L$ be
a field of characteristic $0$.
An orthogonal 
$L$-representation of $\pi$
is a triple
$(V,b,\rho)$
consisting of
an $L$-vector space
$V$ of finite dimension,
a non-degenerate symmetric
bilinear form
$b:V\otimes V\to L$
and a representation
$\rho:\pi\to O(V,b)$
to the orthogonal group
$O(V,b)$.
If $\pi$ is
a profinite group,
we assume that
$L$ is a finite
extension of
the $\ell$-adic field
${\mathbb Q}_\ell$
and that $\rho$
is continuous.

Let $V$ be an orthogonal
representation of $\pi$.
The first Stiefel-Whitney class
$sw_1(V)\in 
H^1(\pi,{\mathbb Z}/2{\mathbb Z})$
is the determinant
$\det \rho\colon
\pi\to \{\pm1\}
\subset L^\times$
regarded as
an element of
$H^1(\pi,{\mathbb Z}/2{\mathbb Z})
=
Hom(\pi,\{\pm1\})$.
The second Stiefel-Whitney class
$sw_2(V)\in H^2(\pi,{\mathbb Z}/2{\mathbb Z})$
is defined as follows.
Let 
\begin{equation}
\begin{CD}
1@>>>{\mathbb Z}/2{\mathbb Z}@>>>
\widetilde {\bf O}(V)@>>>
{\bf O}(V)@>>>1
\end{CD}\label{eqO}
\end{equation}
be the central extension
of the algebraic group
${\bf O}(V)={\bf O}(V,b)$ by ${\mathbb Z}/2{\mathbb Z}$
defined by using the
Clifford algebra ${\rm Cl}(V)$
as in \cite{Fr,EKV,S2}.
The canonical map $\widetilde {\bf O}(V)\to
{\bf O}(V)$ sends $x\in \widetilde {\bf O}(V)
\subset {\rm Cl}(V)$
to the automorphism $v\mapsto I(x)vx^{-1}$
where $I$ denotes the automorphism
of the Clifford algebra defined by
the multiplication by $-1$ on the odd part.
By pulling back the central extension
(\ref{eqO}) by $\rho$,
we obtain a central extension
of $\pi$ by ${\mathbb Z}/2{\mathbb Z}$.
The Stiefel-Whitney class
$sw_2(V)\in H^2(\pi,{\mathbb Z}/2{\mathbb Z})$
is defined as the class of
this central extension
\cite[{\sc Chap}.\ XIV Theorem 4.2]{CE}.

For the orthogonal sum
of orthogonal representations,
we have
$sw_1(V\oplus V')=
sw_1(V)+sw_1(V')$
and 
$sw_2(V\oplus V')=
sw_2(V)+
sw_1(V)sw_1(V')+
sw_2(V')$.
If we introduce the notation
$sw(V)
=1+sw_1(V)+sw_2(V)$,
the equalities
are rewritten as
$sw(V\oplus V')=
sw(V)\cdot sw(V')$
\cite[Lemma 2.1]{S2}.

If an orthogonal
representation $V$
of $\pi$
admits a direct
sum decomposition
$V=W\oplus W'$
by $\pi$-stable
and isotropic subspaces,
the Stiefel-Whitney
classes are
computed as follows.
For a character
$\chi\colon \pi
\to L^\times$,
we define
$\bar c_1(\chi)
\in
H^2(\pi,{\mathbb Z}/
2{\mathbb Z})$
to be
the class of
the pull-back
by $\chi$
of the central extension
\begin{equation}
\begin{CD}
1@>>> {\mathbb Z}/
2{\mathbb Z}
@>>> {\mathbb G}_m
@>2>>
{\mathbb G}_m
@>>> 1.
\end{CD}\label{eqO1}
\end{equation}
For a finite dimensional
$L$-representation
$V$ of $\pi$,
we put
$\bar c_1(V)
=
\bar c_1(\det V)$.

\begin{lm}\label{lmc1}
Let $W$ be
an $L$-vector space
of finite dimension
and we define
a symmetric non-degenerate
bilinear form on
the direct sum
$V=W\oplus W^\vee$
with the dual
by the canonical pairing.
We regard
$GL(W)$ as a subgroup
of $O(V)$
by $g\mapsto g\oplus g^{\vee -1}$.
Then, 
the pull-back
of {\rm (\ref{eqO})}
by the inclusion
$GL(W)\to O(V)$
is isomorphic
to the pull-back
of {\rm (\ref{eqO1})}
by 
$\det\colon
GL(W)\to {\mathbb G}_m$.
\end{lm}

\begin{proof}
Since the algebraic
group $SL(W)$
is connected
and simply connected,
the pull-back
of {\rm (\ref{eqO})}
by the inclusion
$GL(W)\to O(V)$
is isomorphic
to the pull-back
of a central extension of
${\mathbb G}_m$
by 
$\det\colon
GL(W)\to {\mathbb G}_m$.
We define
a section 
${\mathbb G}_m
\to GL(W)$
by taking a line in $W$.
Then, it is reduced
to the case
where $W$ is a line.

We assume $\dim W=1$
and identify
$GL(W)={\mathbb G}_m$.
Let $e\in W$ be a basis
and $f\in W^\vee$
be the linear form
defined by $f(e)=1/2$.
Then, the map
${\mathbb G}_m
\to \widetilde
{\bf O}(V)$
defined by
$a\mapsto
a+(1/a-a) f\cdot e$
is a group homomorphism
since $f\cdot e$ is
an idempotent.
Further,
it makes the diagram
$$\begin{CD}
{\mathbb G}_m
@>2>>
{\mathbb G}_m&
=GL(W)\\
@VVV &@VVV\\
\widetilde
{\bf O}(V)
@>>>&
{\bf O}(V)
\end{CD}$$
commutative
since $e\cdot e=0$ and
$f\cdot e\cdot f=f$.
Hence the assertion
follows.
\end{proof}

\begin{cor}\label{corc1}
Let $V$ be
an orthogonal 
$L$-representation
of $\pi$
and $V=W\oplus W'$ be
a decomposition
by $\pi$-stable isotropic
subspaces.
Then,
we have
$sw_2(V)=\bar c_1(W)$
in $H^2(\pi,
{\mathbb Z}/
2{\mathbb Z})$.
\end{cor}

\begin{proof}
By the assumption,
we may identify
$W'$ with the dual
space of $W$.
Then,
it follows
from Lemma \ref{lmc1}.
\end{proof}

We compute the Stiefel-Whitney
class of the twist
by a character of order $2$.
We prepare lemmas on
central extensions.

\begin{lm}\label{lmcup}
Let $n>1$ be an integer.
Let $1\to {\mathbb Z}/
n{\mathbb Z}\to
\widetilde G
\to G\to 1$
and $1\to {\mathbb Z}/
n{\mathbb Z}\to
\widetilde G'
\to G'\to 1$
be central extensions
and $\chi\colon
G\to {\mathbb Z}/
n{\mathbb Z}$
and $\chi'\colon
G'\to {\mathbb Z}/
n{\mathbb Z}$
be characters.
We define a new group
structure
on the quotient
\setcounter{equation}0
\begin{equation}
E=
(\widetilde G\times \widetilde G')/
{\rm Ker}(+\colon
{\mathbb Z}/
n{\mathbb Z}\oplus
{\mathbb Z}/
n{\mathbb Z}
\to {\mathbb Z}/
n{\mathbb Z})
\label{eqUn}
\end{equation}
by $\overline{(g,g')}
\cdot
\overline{(h,h')}
=
(\chi(h)\cdot\chi'(g'))
\overline{(gh,g'h')}$
where $\cdot$
in the right hand side
denotes the multiplication
of the ring ${\mathbb Z}/
n{\mathbb Z}$.

Then the class
$[E]\in H^2(G\times G',
{\mathbb Z}/
n{\mathbb Z})$
of the central extension
$1\to {\mathbb Z}/
n{\mathbb Z}
\to E\to 
G\times G'\to 1$
is equal to the sum
${\rm pr}_1^*
[\widetilde G]+
{\rm pr}_2^*
[\widetilde G']+
{\rm pr}_1^*
\chi\cup
{\rm pr}_2^*
\chi'$.
\end{lm}

\begin{proof}
Let $U({\mathbb Z}/
n{\mathbb Z})\subset
GL_3({\mathbb Z}/
n{\mathbb Z})$ 
be the subgroup 
consisting of unipotent
upper triangular matrices.
If we put 
$G=G'={\mathbb Z}/
n{\mathbb Z},
\widetilde G=
\widetilde G'=
{\mathbb Z}/
n{\mathbb Z}\oplus
{\mathbb Z}/
n{\mathbb Z}$
and $\chi=\chi'={\rm id}$,
then 
$U({\mathbb Z}/
n{\mathbb Z})$
is obtained 
as $E$
defined in (\ref{eqUn}).
It is a central extension
of 
${\mathbb Z}/
n{\mathbb Z}
\times
{\mathbb Z}/
n{\mathbb Z}$
by 
${\mathbb Z}/
n{\mathbb Z}$
and its class
in $H^2(
{\mathbb Z}/
n{\mathbb Z}
\times
{\mathbb Z}/
n{\mathbb Z},
{\mathbb Z}/
n{\mathbb Z})$
is
${\rm pr}_1
\cup
{\rm pr}_2$.

The class of 
the central extension $E_0$
with the same underlying
set as $E$ and 
with the group structure
without modification
is equal to the sum
${\rm pr}_1^*
[\widetilde G]+
{\rm pr}_2^*
[\widetilde G']$.
Since $E$ is obtained
by modifying $E_0$
by the pull-back of
$U({\mathbb Z}/
n{\mathbb Z})$
by $\chi\times \chi'$,
the assertion follows.
\end{proof}

\begin{lm}\label{lmt2}
Let $b$ be a non-degenerate
symmetric bilinear form
on an $L$-vector space $V$
of finite dimension $n\ge 1$.

{\rm 1.}
Let $1\to \mu_2
\to A\to \mu_2\to 1$ 
be the central extension
defined as the pull-back
of {\rm (\ref{eqO})}
by the inclusion
$\mu_2\to {\mathbf O}(V)$.
Then, the class $[A]
\in H^2(\mu_2,{\mathbb Z}/2
{\mathbb Z})\simeq
{\mathbb Z}/2
{\mathbb Z}$
is $\displaystyle{
\binom n2}$.

{\rm 2.}
We apply the construction
{\rm (\ref{eqUn})}
to the central extensions
$1\to \mu_2
\to A\to \mu_2\to 1$ 
and $1\to \mu_2
\to \widetilde {\mathbf O}(V)
\to {\mathbf O}(V)\to 1$ 
and the characters
${\rm id}^{n-1}\colon
\mu_2\to \mu_2$
and
$\det\colon
{\mathbf O}(V)\to \mu_2$
to define a central extension
$E$ of
${\mathbf O}(V)\times \mu_2$
by $\mu_2$.
Then, the inclusions induce
a morphism
$E\to
\widetilde {\mathbf O}(V)$
of groups.
\end{lm}

\begin{proof}
After extending
$L$ if necessary,
we take an orthonormal basis
$x_1,\ldots,x_n$ of $V$.

1. 
Since the subgroup 
$A$ is generated by
$x_1\cdots x_n$ and $-1$,
the assertion follows
from
$(x_1\cdots x_n)^2=
(-1)^{\binom n2}$.

2.
For $x\in V$,
we have an equality
$x_1\cdots x_n\cdot x
=(-1)^{n-1}x\cdot x_1\cdots x_n$
in the Clifford algebra ${\rm Cl}(V)$.
Since the algebraic group
$\widetilde {\mathbf O}(V)$
is generated by
$V\cap \widetilde {\mathbf O}(V)
\subset {\rm Cl}(V)$,
the assertion follows.
\end{proof}

\begin{cor}\label{cort2}
Let $\rho\colon
\pi\to O(V)$
be an orthogonal representation
of degree $n$
and $\chi\colon
\pi\to \mu_2$ be
a character of order $2$.
We regard
$\det \rho$ and $\chi$
as elements in $H^1(\pi,
{\mathbb Z}/
2{\mathbb Z})$.
Then, we have
\setcounter{equation}0
\begin{equation}
sw_2(\rho\otimes \chi)=
sw_2(\rho)+
(n-1)\det \rho\cup \chi
+\binom n2 \chi\cup \chi
\label{eqt2}
\end{equation}
in $H^2(\pi,
{\mathbb Z}/
2{\mathbb Z})$
and
$\det (\rho\otimes \chi)=
\det \rho +n\chi$
in $H^1(\pi,
{\mathbb Z}/
2{\mathbb Z})$.
\end{cor}

\begin{proof}
By Lemma \ref{lmt2}.2,
the Stiefel-Whitney class
$sw_2(\rho\otimes \chi)$
is the class of the
central extension defined
as the pull-back of $E$
in Lemma \ref{lmt2}.2 by
$\rho\times \chi
\colon \pi \to
{\mathbf O}(V)
\times \mu_2$.
Hence the equation 
(\ref{eqt2}) follows
from Lemmas \ref{lmcup} and
\ref{lmt2}.1.
The assertion for det is clear.
\end{proof}

We generalize the definitions
to graded case.
Let $V^\bullet
=\bigoplus_{q\in 
{\mathbb Z}}V^q$
be a graded 
$L$-representation.
For each integer
$q\in {\mathbb Z}$,
we assume $V^q$
is a finite dimensional
$L$-vector space
equipped with
a continuous representation
of $\pi$
and $V^q=0$
except for finitely many $q$.
We assume
that $V^0$
is an orthogonal representation
and that,
for each $q>0$,
$V^q\oplus V^{-q}$
is equipped with
a $\pi$-invariant
$(-1)^q$-symmetric
non-degenerate form
such that
$V^q$ and $V^{-q}$
are totally isotropic.
Then, we define
$sw_1(V^\bullet)
\in
H^1(\pi,
{\mathbb Z}/
2{\mathbb Z})$ to
be $\det(V^0)$
and
$sw_2(V^\bullet)
\in
H^2(\pi,
{\mathbb Z}/
2{\mathbb Z})$ by
\begin{equation}
sw_2(V^\bullet)
=sw_2(V^0)
+\sum_{q<0}
\bar c_1(V^q).
\label{eqdfswg}
\end{equation}
The definition
is equivalent to
the equality
$$1+sw_1(V^\bullet)
+sw_2(V^\bullet)
=
(1+\det V^0
+sw_2(V^0))
\cdot
\prod_{q<0}
(1+\bar c_1(V^q))^{(-1)^q}
$$
in $1+
H^1(\pi,
{\mathbb Z}/
2{\mathbb Z})+
H^2(\pi,
{\mathbb Z}/
2{\mathbb Z})$.

If $K$ is a field and $\pi$ 
is the absolute Galois group
$G_K=\Gal (\bar K/K)$,
the Stiefel-Whitney classes
$sw_i(V)$ are defined in
the Galois cohomology
$H^i(K,{\mathbb Z}/2{\mathbb Z})
=H^i(G_K,{\mathbb Z}/2{\mathbb Z})$
for $i=1,2$.
If the characteristic of
$K$ is not $2$,
we identify 
$H^1(K,{\mathbb Z}/2{\mathbb Z})$
with $K^\times/(K^\times)^2$.
For $a\in K^\times$,
let $\{a\}\in H^1(K,{\mathbb Z}/2{\mathbb Z})$
denote its class.
For $a,a'\in K^\times$,
let $\{a,a'\}\in H^2(K,{\mathbb Z}/2{\mathbb Z})$ 
denote
the cup-product 
$\{a\}\cup \{a'\}$.

\begin{lm}\label{lmcchi}
{\rm 1.}
Let $\ell$ be a prime number
and
$\chi_\ell\colon
\pi_1({\rm Spec}\
{\mathbb Z}[\frac 1\ell])^{\rm ab}
\to {\mathbb Q}_\ell^\times$
be the $\ell$-adic cyclotomic
character of
the abelianized
algebraic fundamental group.
Then, 
$c_\ell=
\bar c_1(\chi_\ell)$
is the generator
of the cyclic group
$H^2(\pi_1({\rm Spec}\
{\mathbb Z}[\frac 1\ell])^{\rm ab},
{\mathbb Z}/2
{\mathbb Z})$ of order $2$.

{\rm 2.}
Let $K$ be a field of
characteristic $\neq 2$
and $\chi\colon
G_K\to \{\pm1\}
\subset L^\times$
be a character.
Then, 
we have
$\bar c_1(\chi)
=\chi\cup\{-1\}.$
\end{lm} 

\begin{proof}
{\rm 1.}
Since
$\chi_\ell$
defines
an isomorphism
$\pi_1({\rm Spec}\
{\mathbb Z}[\frac 1\ell])^{\rm ab}
\to {\mathbb Z}_\ell^\times$,
the profinite group
$\pi_1({\rm Spec}\
{\mathbb Z}[\frac 1\ell])^{\rm ab}$
is isomorphic
to the product of
${\mathbb Z}_\ell$
with a cyclic group $C$ of
even order.
Hence, 
$H^2(\pi_1({\rm Spec}\
{\mathbb Z}[\frac 1\ell])^{\rm ab},
{\mathbb Z}/2
{\mathbb Z})$ is
of order $2$
and is generated
by the pull-back of
the class of the unique
central extension of
$C$ by ${\mathbb Z}/2
{\mathbb Z}$.

{\rm 2.}
Since
$\bar c_1({\rm id}),
{\rm id}
\cup
{\rm id}\in
H^2({\mathbb Z}/
2{\mathbb Z},
{\mathbb Z}/
2{\mathbb Z})$
are the unique non-trivial
element,
we have
$\bar c_1({\rm id})
={\rm id}
\cup
{\rm id}$.
Hence,
we obtain
$\bar c_1(\chi)
=
\chi\cup\chi$.
Since
$\{a,a\}=
\{a,-1\}$
for $a\in K^\times$
corresponding to $\chi
\in H^1(K,,
{\mathbb Z}/
2{\mathbb Z})$,
the assertion
follows.
\end{proof}

For a field $K$ of
characteristic different from $\ell$,
let 
\setcounter{equation}0
\begin{equation}
c_\ell
=
\bar c_1(\chi_\ell)
\in H^2(K,
{\mathbb Z}/
2{\mathbb Z})
\label{eqcl}
\end{equation} 
also denote
the pull-back
by the canonical map
$G_K\to
\pi_1({\rm Spec}\
{\mathbb Z}[\frac 1\ell])^{\rm ab}$.
If $K={\mathbb Q}$,
then $c_\ell\in
H^2({\mathbb Q},
{\mathbb Z}/2{\mathbb Z})$
is the unique element
ramifying exactly at $\ell$ and $\infty$.
We have $c_2=
\{-1,-1\}$ for example.
If $K$ is of positive
characteristic,
we have $c_\ell=0$.

Let $K$ be a field
of characteristic different
from 2.
We call a pair $(D,b)$
of a $K$-vector space $D$
of finite dimension
and a non-degenerate symmetric
bilinear form $b:D\otimes_KD\to K$
a quadratic $K$-vector space.
The first
and the second
Hasse-Witt classes
$hw_1(D,b)\in H^1(K,{\mathbb Z}/2{\mathbb Z})$
and
$hw_2(D,b)\in H^2(K,{\mathbb Z}/2{\mathbb Z})$
are defined as follows \cite{M}.
Let 
$x_1,\ldots,x_r$ 
be an orthogonal 
basis of $D$
and put $a_i=b(x_i,x_i)$.
The discriminant
of $D$ is defined by
$\disc D
=\prod_{i=1}^r
a_i\in K^\times/
K^{\times 2}$.
We define
the classes
by
$$hw_1(D)
=
\{\disc D\}=
\sum_{i=1}^r\{a_i\}
\quad\text{and}
\quad
hw_2(D)
=
\sum_{1\le i<j\le r}\{a_i,a_j\}.$$
They are independent of
the choice of
an orthogonal basis.

For the orthogonal sum
of quadratic
vector spaces,
we have
$hw_1(D\oplus D')=
hw_1(D)+hw_1(D')$
and 
$hw_2(D\oplus D')=
hw_2(D)+
hw_1(D)hw_1(D')+
hw_2(D')$.
If we introduce the notation
$hw(D)
=1+hw_1(D)+hw_2(D)$,
the equalities
are rewritten as
$hw(D\oplus D')=
hw(D)\cdot hw(D')$.

\begin{lm}\label{lmitp}
Let $(D,b)$ be a quadratic $K$-vector space.

{\rm 1}.
Let
$D^-$ be a totally isotropic subspace
of dimension $r$.
Then the restriction of $b$
to $D^+=(D^-)^\perp$
induces a non-degenerate symmetric bilinear
form $b^0$ on $D^0=D^+/D^-$
and we have
$$
hw_2(D)=
hw_2(D^0)+
r\{-1,\disc D^0\}+
\binom r2\{-1,-1\}.$$

{\rm 2.}
For $a\in K^\times$,
we have
\begin{align*}
hw(D,a\cdot b)=
1+&\ {\dim D}\cdot \{a\}+
\binom {\dim D}2
\cdot\{a,a\}\\
&+\disc D
+
({\dim D}-1)\{a,
\disc D\}
+
hw_2(D,b).
\end{align*}
\end{lm}

\begin{proof}
{\rm 1.}
The quadratic space
$D$ is isomorphic
to the orthogonal direct sum
of $D^0$ with $r$-copies
of the hyperbolic plane.
Hence we have
$hw(D)=
hw(D^0)(1+\{-1\})^r$.

{\rm 2.}
Clear from the definition.
\end{proof}

We generalize the definitions
to graded case.
Let $D^\bullet
=\bigoplus_{q\in 
{\mathbb Z}}D^q$
be a graded 
$K$-vector space.
For each integer
$q\in {\mathbb Z}$,
we assume $D^q$
is a finite dimensional
$K$-vector space
and $D^q=0$
except for finitely many $q$.
We assume
that $D^0$
is a quadratic vector space
and that,
for each $q>0$,
$D^q\oplus D^{-q}$
is equipped with
a $(-1)^q$-symmetric
non-degenerate form
such that
$D^q$ and $D^{-q}$
are totally isotropic.
We put
$r=\sum_{q<0}(-1)^q\dim D^q$.
Then, we define
$hw_1(D^\bullet)
\in
H^1(K,{\mathbb Z}/
2{\mathbb Z})$ 
and
$hw_2(D^\bullet)
\in
H^2(K,{\mathbb Z}/
2{\mathbb Z})$ by
\setcounter{equation}0
\begin{align}
hw_1(D^\bullet)
=&\
\disc (D^0)
+r\cdot\{-1\}
\label{eqdfhwg1}
\\
hw_2(D^\bullet)
=&\
hw_2(D^0)
+
r\cdot
\{-1,\disc D\}
+
\binom
r2
\cdot
\{-1,-1\}.
\label{eqdfhwg2}
\end{align}
The definition
is equivalent to
the equality
$$1+hw_1(D^\bullet)
+hw_2(D^\bullet)
=
(1+hw_1(D^0)+hw_2(D^0))
\cdot
\prod_{q<0}
(1+\{-1\})^{(-1)^q\dim D^q}
$$
in $1+
H^1(K,{\mathbb Z}/
2{\mathbb Z})+
H^2(K,{\mathbb Z}/
2{\mathbb Z})$.

\section{Conjecture}\label{scn}

To formulate
the conjecture,
we prove a preliminary
result on the
determinant
of $\ell$-adic cohomology.

\begin{lm}\label{lm1}
Let $S$ be a connected
normal scheme, $f\colon
X\to S$ be a proper smooth
morphism
and $q\ge0$
be an integer.
For a prime number 
$\ell$ invertible on $S$,
let $\chi_\ell:
\pi_1(S)^{\rm ab}
\to {\mathbb Q}_\ell^\times$
denote the $\ell$-adic
cyclotomic character.

{\rm 1.} {\rm (\cite[Corollary 2.2.3]{Su})}
The rank $b_{\acute et, q}$ of
the smooth
${\mathbb Q}_\ell$-sheaf
$R^qf_*{\mathbb Q}_\ell$
is independent of 
$\ell$ invertible on $S$.
If $q$ is odd,
then the Betti number $b_{\acute et, q}$ is even.

{\rm 2.} {\rm (\cite[Corollary 3.3.5]{Su})}
The character $e_q:
\pi_1(S)^{\rm ab}
\to
{\mathbb Q}_\ell^\times$
defined by
\setcounter{equation}0
\begin{equation}
e_q=
\det R^qf_*
{\mathbb Q}_\ell
\cdot
\chi_\ell^{
q\cdot b_{\acute et, q}/2}
\label{eqeq}
\end{equation}
is independent of
$\ell$ invertible on $S$
and takes values in $\{\pm1\}$.
Further if $q$ is odd,
the character $e_q$ is trivial.
\end{lm}

\begin{proof}
1. By a standard limit argument,
we may assume $S$ is of finite type
over ${\mathbb Z}$. Then, it follows from the
proper base change theorem,
the Weil conjecture and
\cite[Corollary 2.2.3]{Su}.

2.
First, we consider the case
where $S={\rm Spec}\ k$
for a finite field $k$.
Then by the Weil conjecture
proved by Deligne,
$\det(Fr_k\colon
H^q(X_{\bar k},{\mathbb Q}_\ell))$
is a rational integer independent of $\ell$
and is of absolute value
$qb_{\acute et, q}/2$-th power of
${\rm Card}\ k$.
Hence 
we have
$\det(Fr_k\colon
H^q(X_{\bar k},{\mathbb Q}_\ell))=
\pm ({\rm Card}\ k)^{qb_{\acute et, q}/2}$.
Further by \cite[Corollary 3.3.5]{Su}),
the sign is positive if $q$ is odd.
Thus the assertion follows in this case.

We prove the general case.
By replacing $S$ by a dense open,
we may assume $S$ is affine.
By a standard limit argument,
we may assume $S$ is of finite type over
${\mathbb Z}$.
By the Chebotarev density theorem 
\cite[Theorem 7]{ZL},
\cite[Theorem 9.11]{STai},
the reciprocity map
$\bigoplus_{s\in S_0}
{\mathbb Z}\to 
\pi_1(S)^{\rm ab}$
has dense image,
where $S_0$ denotes
the set of closed points in $S$.
Since it is proved for the spectrum of a finite field,
it follows that
the character 
$\det R^qf_*
{\mathbb Q}_\ell$
is independent of $\ell$ and 
its square is equal to 
$\chi_\ell^{-
q\cdot b_{\acute et, q}}$.
Further if $q$ is odd,
the character 
$\det R^qf_*
{\mathbb Q}_\ell$
itself is equal to 
$\chi_\ell^{-
q\cdot b_{\acute et, q}/2}$.
\end{proof}

Let $X$ be a proper
smooth scheme 
of even dimension $n$ over 
a field $K$.
Let $\ell$ be a prime number
different from the characteristic of $K$.
The cup-product 
defines a non-degenerate 
symmetric bilinear form
on $V=H^n(X_{\bar K},
{\mathbb Q}_\ell)(\frac n2)$.
We consider $V$
as an orthogonal representation
of the absolute Galois group
$G_K=\Gal(\bar K/K)$ and
define 
\setcounter{equation}0
\addtocounter{thm}1
\begin{equation}
sw_2(H^n_\ell(X))
\in H^2(K,{\mathbb Z}/2{\mathbb Z})
\label{eqsw}
\end{equation}
to be its second Stiefel-Whitney class.

Assume the characteristic
of $K$ is not 2.
The cup-product
defines a non-degenerate 
symmetric bilinear form
also on $D=H^n_{dR}(X/K)$.
We consider $D$
as a quadratic $K$-vector space
and define 
\begin{equation}
hw_2(H^n_{dR}(X))
\in H^2(K,{\mathbb Z}/2{\mathbb Z})
\label{eqhw}
\end{equation}
to be its second Hasse-Witt class.
We also put
\begin{equation}
d_X=\disc 
H^n_{dR}(X/K)
=hw_1(H^n_{dR}(X))
\in H^1(K,{\mathbb Z}/2{\mathbb Z}).
\label{eqhw1}
\end{equation}

To state a conjecture on
the relation between
$sw_2(H^n_\ell(X))$
and $hw_2(H^n_{dR}(X))$,
we introduce auxiliary invariants.
For an integer $q$,
we put 
\begin{equation}
b_{\acute et,q}=
\dim H^q(X_{\bar K},{\mathbb Q}_\ell)
\quad
\text{and} 
\quad
b_{dR,q}=
\dim H^q_{dR}(X/K).
\label{eqbq}
\end{equation}
If $K$ is of characteristic 0,
we have $b_{\acute et,q}=
b_{dR,q}$ and we simply write
them $b_q$. 

\begin{cn}\label{cn}
Let $X$ be a proper and smooth 
scheme of even dimension $n$
over a field $K$ of
characteristic $\neq 2,\ell$.
For an integer $q\ge 0$,
we regard the character
$e_q:G_K\to \{\pm1\}$
{\rm (\ref{eqeq})}
as an element of
$Hom(G_K,\{\pm1\})=
H^1(K,{\mathbb Z}/2{\mathbb Z})$ and put
$e=\sum_{q<n} e_q$,
\setcounter{equation}0
\begin{equation}
r=\sum_{q<n}(-1)^q b_{dR,q},
\qquad
\eta=
\sum_{q<\frac n2}
(-1)^q\left(\frac n2-q\right)
\chi(X,\Omega^q_{X/K}).
\label{eqr}
\end{equation}
If the characteristic of $K$ is $0$,
we further define
\begin{equation}
\beta=
\frac12\sum_{q<n}
(-1)^q(n-q)b_q.
\label{eqbeta}
\end{equation}

Then we have an equality
\begin{align}
&\ sw_2(H^n_\ell(X))+ \{e,-1\}
+\beta\cdot c_\ell
\label{eqcn}
\\
=&\ hw_2(H^n_{dR}(X))
\nonumber \\
&\qquad +
\begin{cases}\displaystyle{
r\{d_X,-1\}+\binom r2\{-1,-1\}}
& \text{if }n\equiv 0 \bmod 4,
\nonumber\\\displaystyle{
(r+b_{dR,n}-1)\{d_X,-1\}+
\binom {r+b_{dR,n}}2\{-1,-1\}}
& \text{if }n\equiv 2 \bmod 4
\end{cases}
\nonumber\\
&\qquad + \{2,d_X\}+
\eta\cdot (c_\ell-c_2)
\nonumber
\end{align}
in $H^2(K,{\mathbb Z}/2{\mathbb Z})$.
\end{cn}

An apparently simpler reformulation of the conjecture
will be given at Corollary \ref{corg}.
A generalization for a family is stated in
Conjecture \ref{cn2}
after defining the Stiefel-Whitney class
for a symmetric perfect complex
in Section \ref{sswcx}.

\begin{rmk}
Let $X$ be a projective smooth scheme
of even dimension $n$
over a field $K$ of
characteristic $\neq 2,\ell$.
We regard the character
$e_n\colon
G_K\to \{\pm 1\}$
defined by
$e_n(\sigma)=
\det(\sigma:H^n(X_{\bar K},
{\mathbb Q}_\ell)(\frac n2))$
as an element
in $H^1(K,{\mathbb Z}/2{\mathbb Z})
=Hom(G_K,\{\pm 1\})$.
Then the equality
\setcounter{equation}0
\begin{equation}
e_n=
\{d_X\}+
\begin{cases}
r\cdot\{-1\}
&\quad\text{ if } n\equiv0 \bmod 4,\\
(r+b_{dR,n})\cdot\{-1\}
&\quad\text{ if } n\equiv2 \bmod 4
\end{cases}
\label{eqdet}
\end{equation}
in $H^1(K,{\mathbb Z}/2{\mathbb Z})$
is proved
in {\rm \cite[Theorem 2]{S1}}.
Conjecture {\rm \ref{cn}}
is a degree $2$ version
of the equality 
{\rm (\ref{eqdet})}.
\end{rmk}

In this paper,
we prove the following
evidence for Conjecture \ref{cn}.

\begin{thm}\label{thm1}
Let $X$ be a proper and smooth 
scheme of even dimension $n$
over a field $K$ of
characteristic $\neq 2,\ell$.
Conjecture {\rm \ref{cn}}
is true
in the following cases.

{\rm 1.} 
$K$ is a finite extension of 
${\mathbb Q}_p$, $p\neq 2,\ell$
and there exists a 
projective regular
flat model $X_{{\mathcal  O}_K}$ over the
integer ring ${\mathcal  O}_K$
such that
the closed fiber has
at most ordinary double points
as singularities.

{\rm 2.} 
$K$ is a finite extension of 
${\mathbb Q}_p$, $p=\ell> n+1$
and there exists a proper smooth 
model $X_{{\mathcal  O}_K}$ over the
integer ring ${\mathcal  O}_K$.

{\rm 3.} 
$K={\mathbb R}$ and
$X$ is projective.

{\rm 4.} 
$K\supset \bar {\mathbb Q}$.

{\rm 5.} 
$X$ is a smooth
hypersurface in ${\mathbb P}^{n+1}_K$
and $l> n+1$.
\end{thm}

Theorem \ref{thm1} implies that 
the Stiefel-Whitney class
$sw_2(H^n_\ell(X))$
may depend on $\ell$.

\begin{ex}\label{ex}
Let $p$ be a prime
and $X$ be a proper
smooth variety of even dimension $n$
over $K={\mathbb Q}_p$ with good reduction.
Then for any prime $\ell\neq p$,
we have
$sw_2(H^n_\ell(X))=
hw_2(H^n_{dR}(X))=0$
and $d,e\in H^1({\mathbb Q}_p,
{\mathbb Z}/2{\mathbb Z})$
are unramified.
Further if $p>3$ and $A$ is
an abelian surface,
Theorem
{\rm \ref{thm1}.2}
implies 
$sw_2(H^2_p(A))
=c_p \neq0$
since
$\beta=\frac12(2-4)\equiv 1\bmod 2$
and $\eta=0$.

Similarly, for an abelian 
surface $A$ over ${\mathbb R}$,
we have
$$
sw_2(H^2_\ell(A))\neq
hw_2(H^2_{dR}(A))=0.$$
\end{ex}

If $n=0$,
the assertion 5 in Theorem \ref{thm1}
is nothing but
the following theorem of Serre
since $e=1$
and $r,\beta, \eta=0$ in this case.
The proof of the assertion 5
gives a new proof of
the formula of Serre.

\begin{thm}
[{{\rm \cite[Theorem 1]{Se}}}]
\label{thm2}
Let $K$ be a field
of characteristic $\neq2$
and $L$ be a finite separable
extension of $K$.
We consider 
$V={\mathbb Q}^{{\rm Mor}_K(L,\bar K)}$
as an orthogonal representation 
of $G_K$
and $D=L$
as a quadratic $K$-vector space
with a non-degenerate
symmetric bilinear form
$(x,y)\mapsto
\text{\rm Tr}_{L/K}(xy)$.
Let $d_{L/K}\in H^1(K,
{\mathbb Z}/2{\mathbb Z})$ 
be the discriminant
of $D=L$.
Then we have
$$sw_2({\mathbb Q}^{{\rm Mor}_K(L,\bar K)})=
hw_2(L,\text{\rm Tr}_{L/K}(xy))
+\{2,d_{L/K}\}.$$
\end{thm}

\begin{proof}
We put $X={\rm Spec}\ L$.
Then the $\ell$-adic representation
$H^0(X_{\bar K},{\mathbb Q}_\ell)$
is the extension of scalars of
the continuous representation
${\mathbb Q}^{{\rm Mor}_K(L,\bar K)}$.
Hence, the left hand side
is equal to
$sw_2(H^0_\ell(X))$.
Since the cup-product on
$L=H^0_{dR}(X/K)$
is the multiplication of $L$
and the trace map
${\rm Tr}\colon H^0_{dR}(X/K)=L\to K$
is the usual trace map
for the extension $L$ over $K$,
the first term of the right hand
side equals
$hw_2(H^0_{dR}(X))$.
Hence the assertion follows from the case where $n=0$
of the assertion 5 in Theorem \ref{thm1}.
\end{proof}

We give some
simplified versions
of the formula (\ref{eqcn}).

\begin{lm}\label{lm'}
Let $hw_2(H^n_{dR}(X)')$
denote the second 
Hasse-Witt class
of $H^n_{dR}(X)$
with the symmetric bilinear
form multiplied by
$(-1)^{n/2}$
and
$d'_X=
hw_1(H^n_{dR}(X)')$
denote the discriminant.
Then, if $n\equiv 2\bmod 4$,
we have
\setcounter{equation}0
\begin{align}
hw_2(H^n_{dR}(X)')
&= 
hw_2(H^n_{dR}(X))
+
(b_{dR,n}-1)\cdot \{d_X,-1\}+
\binom {b_{dR,n}}2\{-1,-1\},
\label{eqhw'}\\
d'_X
&=
d_X
+b_{dR,n}\cdot \{-1\}.
\nonumber
\end{align}
\end{lm}

\begin{proof}
It suffices
to apply Lemma \ref{lmitp}.2.
\end{proof}

\begin{cor}\label{cor'}
The equality 
{\rm (\ref{eqcn})} 
is equivalent to
the following:
\setcounter{equation}0
\begin{align}
&\ sw_2(H^n_\ell(X))+ \{e,-1\}
+\beta\cdot c_\ell
\label{eqcn'}
\\
=&\ hw_2(H^n_{dR}(X)')
+
r\{d'_X,-1\}+\binom r2\{-1,-1\}
+ \{2,d_X\}+
\eta\cdot (c_\ell-c_2)
\nonumber
\end{align}
in $H^2(K,{\mathbb Z}/2{\mathbb Z})$.
\end{cor}

The equality 
{\rm (\ref{eqcn})}
is further simplified
if we use
generalizations of
the classes
for graded
orthogonal representations
and graded quadratic forms.
Following the definition
{\rm (\ref{eqdfswg})},
we put
$H^q_\ell(X)=
H^q(X_{\bar K},
{\mathbb Q}_\ell)(\frac n2)$
and
\begin{equation}
sw_2(H^\bullet_\ell(X))
=
sw_2(H^n_\ell(X))
+
\sum_{0\le q<n}
\bar c_1(H^q_\ell(X)).
\label{eqcnsg}
\end{equation}
Similarly,
following
{\rm (\ref{eqdfhwg2})},
we put
\begin{equation}
hw_2(H^\bullet_{dR}(X))
=
hw_2(H^n_{dR}(X)')
+r\cdot
\{-1,d'_X\}
+
\binom
r2
\cdot
\{-1,-1\}
\label{eqcnhg}
\end{equation}
where
$r=\sum_{q<n}(-1)^q b_{dR,q}$
(\ref{eqr}).

\begin{lm}\label{lmg}
We have the following
equality:
\setcounter{equation}0
\begin{equation}
sw_2(H^\bullet_\ell
(X))=
sw_2(H^n_\ell(X))+ \{e,-1\}
+\beta\cdot c_\ell.
\label{eqswcx}
\end{equation}
\end{lm}

\begin{proof}
By the definition of
$e_q$, we have
$\det H^q_\ell(X)=
e_q\cdot
\chi_\ell^{(n-q)b_{dR,q}/2}$.
Hence, 
we obtain
$\sum_{0\le q<n}
\bar c_1(H^q_\ell(X))
=
\bar c_1(e)
+\beta \bar c_1(\chi_{\ell})$
by the definitions
$e=\sum_{q<n} e_q$
and
$\beta=\displaystyle{
\frac12\sum_{q<n}
(-1)^q(n-q)
b_{\text{\it \'et},q}}$.
Thus, the equality
(\ref{eqswcx})
follows from
Lemma \ref{lmcchi}.
\end{proof}

\begin{cor}\label{corg}
The equality 
{\rm (\ref{eqcn})} 
is equivalent to
the following:
\setcounter{equation}0
\begin{equation}
sw_2(H^\bullet_\ell(X))
= 
hw_2(H^\bullet_{dR}(X))
+\{2,d_X\}+
\eta\cdot (c_\ell-c_2).
\label{eqcncx}
\end{equation}
\end{cor}

We compare the de Rham cohomology
with the Hodge cohomology.
We put $n=2m$
and we regard
$H^m(X,\Omega^m_{X/K})$
as a quadratic vector space
by symmetric bilinear form
defined by
the composition
$$\begin{CD}
H^m(X,\Omega^m_{X/K})
\times
H^m(X,\Omega^m_{X/K})
@>{\cup}>>
H^n(X,\Omega^n_{X/K})
@>{\rm Tr}>>
K.
\end{CD}$$

\begin{lm}\label{lmHdR}
The dimensions
$b_{dR,n}=
H^n_{dR}(X/K)$
and
$h^{m,m}=
H^m(X,\Omega^m_{X/K})$
have the same parity.
We put
$b_{dR,n}=
h^{m,m}+2s$.
Then, we have
\setcounter{equation}0
\begin{align}
hw_2(H^n_{dR}(X/K))
=&\
hw_2(H^m(X,\Omega^m_{X/K}))
\label{eqHdR}
\\
&\ +s\cdot
\{-1,hw_1(H^m(X,\Omega^m_{X/K}))\}
+
\binom
s2
\cdot
\{-1,-1\}
\nonumber
\end{align}
and 
$hw_1(H^n_{dR}(X/K))
=
hw_1(H^m(X,\Omega^m_{X/K}))+
s\cdot \{-1\}.$
\end{lm}

\begin{proof}
The Hodge to de Rham spectral sequence
$E_1^{p,q}=
H^q(X,\Omega^p_{X/K})
\Rightarrow
H^r_{dR}(X/K)$
defines decreasing filtrations $F^\bullet$
on $E_1^{m,m}=
H^m(X,\Omega^m_{X/K})$
and on 
$H^n_{dR}(X/K)$.
We choose numbering of
the filtration so that we have
$F^0/F^1=E_\infty^{m,m}$
and put
$s_1=\dim F^1H^n_{dR}(X/K)$
and
$s_2=\dim F^1H^m(X,\Omega^m_{X/K})$.
Then, since the spectral sequence
is compatible with the Serre duality,
we have $F^0=(F^1)^\perp$.
Hence we have $s=s_1-s_2$ and
$hw(H^n_{dR}(X/K))
=
hw(E_\infty^{m,m})
(1+\{-1\})^{s_1}$ and
$hw(H^m(X,\Omega^m_{X/K}))
=
hw(E_\infty^{m,m})
(1+\{-1\})^{s_2}$
by Lemma \ref{lmitp}.1.
Thus, the assertion follows.
\end{proof}

If the characteristic is 0,
Conjecture \ref{cn}
may be restated
as follows.

\begin{lm}\label{lmch0}
Let $X$ be a proper and smooth 
scheme of even dimension $n$
over a field $K$ of
characteristic $0$.

{\rm 1.} 
We put 
\setcounter{equation}0
\begin{equation}
r'=\begin{cases}
b_0-b_2+b_4-\cdots -b_{n-2}
&\quad\text{ if }n\equiv 0\bmod 4\\
-b_0+b_2-b_4+\cdots -b_{n-2}+b_n
&\quad\text{ if }n\equiv 2\bmod 4
\end{cases}
\label{eqr'}
\end{equation}
and
\begin{equation}
h=\sum_{q<\frac n2}
\left(\frac n2-q\right)
\dim H^{n-q}(X,\Omega^q_{X/K}).
\label{eqh}
\end{equation}
Then in this case,
the equality {\rm (\ref{eqcn})}
in Conjecture {\rm \ref{cn}} is 
equivalent to
\begin{align}
&sw_2(H^n_\ell(X))+ \{e,-1\}
\label{eqp}
\\=&\ hw_2(H^n_{dR}(X))
+\{2,d_X\}+
h\cdot (c_\ell-c_2)
+
(r'-\frac n2)\{d_X,-1\}+
\binom {r'}2
\{-1,-1\}.
\nonumber
\end{align}

{\rm 2.} 
Assume further $X$ is projective.
Let
$H^n(X_{\bar K},{\mathbb Q}_\ell)(\frac n2)=
\bigoplus_{q\le n,{\rm even}}P^q$
be the Lefschetz decomposition
into primitive parts
by an ample invertible sheaf
${\mathcal L}$ and put
$P^+=
\bigoplus_{q< n,q\equiv 0\bmod 4}P^q$
and
$P^-=
\bigoplus_{q< n,q\equiv 2\bmod 4}P^q$.
Then, 
we have a congruence
\begin{equation}
r+2\beta
\equiv
-\begin{cases}
\dim P^-
\quad\text{ if } 
n\equiv 0\bmod 4\\
\dim P^+
\quad\text{ if } 
n\equiv 2\bmod 4,
\end{cases}
\label{eqLefr}
\end{equation}
modulo $4$ and 
an equality
\begin{equation}
e=
\begin{cases}
\det P^-
\quad\text{ if } 
n\equiv 0\bmod 4\\
\det P^+
\quad\text{ if } 
n\equiv 2\bmod 4
\end{cases}
\label{eqLefe}
\end{equation}
in $H^1(K,{\mathbb Z}/
2{\mathbb Z})$.
\end{lm}

\begin{proof}
1. 
Since $c_2=\{-1,-1\}$,
it is sufficient to
show the congruences
\begin{align}
\beta&
\equiv
\eta+h
\label{eqcg1}
\\
r'&\equiv
\begin{cases}
r&\quad \text{ if }
n\equiv 0\bmod 4,\\
r+b_n&\quad \text{ if }
n\equiv 2\bmod 4
\end{cases}
\label{eqcg2}
\\
\binom{r'}2&\equiv
\beta+
\begin{cases}
\binom r2
&\quad \text{ if }
n\equiv 0\bmod 4,\\
\binom {r+b_n}2
&\quad \text{ if }
n\equiv 2\bmod 4
\end{cases}
\label{eqcg3}
\end{align}
modulo $2$.
By the Lefschetz principle,
we may assume $K={\mathbb C}$.

We prove (\ref{eqcg1}).
We put $h^{q,p}
=\dim H^p(X,\Omega^q_X)$.
By the Hodge symmetry
and the Serre duality,
we have
$h^{q,p}
=h^{p,q}$
and
$h^{q,p}
=h^{n-q,n-p}$
respectively.
Since the Hodge to de
Rham spectral sequence
degenerates,
we have
$b_q=\sum_{i+j=q}
h^{i,j}$.
By the definition
$\beta=\frac12
\sum_{q<n}(-1)^q
(n-q)b_q$,
we have
\begin{align*}
2\beta
=&\
\sum_{p+q<n}(-1)^{p+q}(n-(p+q))
h^{p,q}\\
=&\
\sum_{p+q<n}(-1)^{p+q}
\left(\frac n2-p\right)
h^{p,q}
+
\sum_{p+q<n}(-1)^{p+q}
\left(\frac n2-q\right)
h^{p,q}.
\end{align*}
Hence by the Hodge symmetry
and the Serre duality,
we obtain
\begin{align*}
\beta
&\
=
\sum_{p+q<n}(-1)^{p+q}
\left(\frac n2-q\right)
h^{p,q}\\
&\
=
\sum_{p+q<n,q<\frac n2}
(-1)^{p+q}
\left(\frac n2-q\right)
h^{p,q}
-
\sum_{p+q>n,q<\frac n2}
(-1)^{p+q}
\left(\frac n2-q\right)
h^{p,q}
.\end{align*}
Thus, we obtain
$\beta+h\equiv
\sum_{q<\frac n2}
(-1)^{p+q}
\left(\frac n2-q\right)
h^{p,q}
=\eta \bmod 2$.

We prove
(\ref{eqcg2}) and
(\ref{eqcg3}).
By Hodge symmetry,
it follows
that the Betti number
$b_q$ is even
for odd $q$.
Hence by the definition
$r=\sum_{q<n}(-1)^q
b_q$ and the definition
of $\beta$ recalled above,
we have
\begin{align}
r+2\beta=&\
\sum_{q<n}(-1)^q
(n-q+1)b_q
\label{eqbr}
\\
\equiv&\
\sum_{q<n,\text{even}}
(n-q+1)b_q
\equiv
r'-
\begin{cases}
0&\quad \text{ if }
n\equiv 0\bmod 4,\\
b_n&\quad \text{ if }
n\equiv 2\bmod 4
\end{cases}
\nonumber
\end{align}
modulo 4.
Hence the congruence
(\ref{eqcg2}) follows.
Since
$\binom{a+2 b}2
\equiv
\binom a2+ b
\bmod 2$,
the congruence
(\ref{eqcg3}) 
also follows.

2.
By (\ref{eqbr}),
we have
$$r+2\beta\equiv
\begin{cases}
(b_0-b_2)+\cdots+(b_{n-4}-b_{n-2})
&\quad \text{ if }
n\equiv 0\bmod 4,\\
-b_0+(b_2-b_4)+\cdots+(b_{n-4}-b_{n-2})
&\quad \text{ if }
n\equiv 2\bmod 4
\end{cases}
$$
modulo 4.
Since $b_q-b_{q-2}=\dim P^q$
for $2\le q\le n-2$
and $b_0=\dim P^0$,
the congruence
(\ref{eqLefr}) follows.

For odd $q$,
$H^q(X_{\bar K},{\mathbb Q}_\ell)$
carries a non-degenerate
alternating form
by hard Lefschetz
and hence 
we have 
$e_q=1$.
For even $2\le q\le n-2$,
we have
$e_q-e_{q-2}=\det P^q$
and
$e_0=\det P^0$.
Since 
$$e=
\begin{cases}
(e_0-e_2)+\cdots+(e_{n-4}-e_{n-2})
&\quad \text{ if }
n\equiv 0\bmod 4,\\
-e_0+(e_2-e_4)+\cdots+(e_{n-4}-e_{n-2})
&\quad \text{ if }
n\equiv 2\bmod 4,
\end{cases}
$$
the equality
(\ref{eqLefe}) is also proved.
\end{proof}

\section{Degenerations}\label{sdeg}

In this section,
we assume that
$K$ is
a complete discrete
valuation field
with residue field $F$
of characteristic $p\neq 2$.
Let $I
=\mathrm {Gal}(\bar K/K^{\rm ur})
\subset
G_K=\mathrm {Gal}(\bar K/K)$
denote the
inertia subgroup
and let $P$
denote the kernel
of  the canonical surjection
$I\to 
\prod_{p'\neq p}
{\mathbb Z}_{p'}(1)$.
Since $P$ is a pro-$p$ group,
the canonical map
$H^q(I/P,{\mathbb Z}/
2{\mathbb Z})\to
H^q(I,{\mathbb Z}/
2{\mathbb Z})$ is an isomorphism
and they are isomorphic
to ${\mathbb Z}/
2{\mathbb Z}$ for $q=0,1$
and
is 0 for $q>1$.
Since the extension
$1\to I/P\to G_K/P\to G_F\to 1$
splits, we have an exact sequence
\setcounter{equation}0
\begin{equation}
\begin{CD}
0@>>>
H^2(F,{\mathbb Z}/
2{\mathbb Z})
@>>>
H^2(K,{\mathbb Z}/
2{\mathbb Z})
@>{\partial}>>
H^1(F,{\mathbb Z}/
2{\mathbb Z})
@>>>0.
\end{CD}
\label{eqpart}
\end{equation}
In this section,
we prove
that the both sides of
(\ref{eqcn}) in Conjecture \ref{cn}
have the same images by
the map $\partial$
in some cases.
In particular, if
$H^2(F,{\mathbb Z}/
2{\mathbb Z})=0$
for example if the residue field
$F$ is finite,
this will imply
Conjecture \ref{cn}
in these cases.

First, we consider
a consequence of the
following elementary lemma.

\begin{lm}\label{lmp2}
Let $K$ be a complete
discrete valuation field
and assume that the characteristic
$p$ of the residue field $F$
is not $2$.

{\rm 1.}
Let $L$ be a finite extension 
of ${\mathbb Q}_\ell$
and $V$ be an orthogonal 
$L$-representation of
$G_K=\mathrm {Gal}(\bar K/K)$.
If the inverse image
$I'\subset I$ of the pro-$2$ part
of ${\mathbb Z}_2(1)
\subset
\prod_{p'\neq p}
{\mathbb Z}_{p'}(1)$
acts
trivially on $V$,
then we have 
$\partial sw_2(V)=0$.

{\rm 2.}
Let $D$ be a quadratic $K$-vector space.
If $D$ has a 
non-degenerate ${\mathcal  O}_K$-lattice,
then we have $\partial hw_2(D)=0$.
\end{lm}

\begin{proof}
1.
The class $sw_2(V)$
is in the image of
$H^2(G_K/I',{\mathbb Z}/2{\mathbb Z})$.
Since $I/I'$ has no pro-2 part,
the canonical map
$H^2(F,{\mathbb Z}/2{\mathbb Z})
=H^2(G_K/I,{\mathbb Z}/2{\mathbb Z})
\to H^2(G_K/I',{\mathbb Z}/2{\mathbb Z})$
is an isomorphism.
Hence the assertion follows
from the exact sequence
(\ref{eqpart}).

2. 
The class
$hw_2(D)$
lies in the image 
$H^2(F,{\mathbb Z}/2{\mathbb Z})
\subset
H^2(K,{\mathbb Z}/2{\mathbb Z})$
of
$\{{\mathcal  O}_K^\times,{\mathcal  O}_K^\times\}$.
\end{proof}

We say a scheme $X$
over
$S={\rm Spec}\ {\mathcal  O}_K$
is generalized semi-stable
if \'etale locally on $X$,
it is \'etale over
${\rm Spec}\ {\mathcal  O}_K
[T_0,\ldots,T_n]/
(T_0^{m_0}\cdots T_r^{m_r}
-\pi)$
for some integer
$0\le r\le n$,
a prime element $\pi$
of $K$
and integers
$m_0,\ldots,m_r\ge 1$
invertible in ${\mathcal  O}_K$.
A scheme $X$ is semi-stable
if and only if
it is generalized semi-stable
and if the
closed fiber
$X_s$ is reduced.

Let $X$ be a generalized semi-stable
scheme over $S$.
Let $\ell$ be a prime number
invertible on $S$
and $R^q\psi{\mathbb Z}_\ell$
be the sheaf of
nearby cycles.
Then, for each geometric point
$\bar x$ of the geometric
closed fiber $X_{\bar s}$,
the action of 
$I$ on the stalk
$R^q\psi{\mathbb Z}_{\ell,\bar x}$
is tamely ramified
for every $q\ge 0$.
More precisely,
if $m_1,\ldots,m_r$
are the multiplicities
of irreducible components
of the closed fiber
$X_{\bar x}\times_S{\rm Spec}\ F
=
\sum_im_iD_i$
of the 
strict henselization
and if $d_x$
denotes the greatest common divisor,
the inertia $I$
acts on
$R^q\psi{\mathbb Z}_{\ell,\bar x}$
through the quotient
$\mu_{d_x}$
\cite[Proposition 6]{St}.

The sheaf 
$\Omega^1_{X/S}(\log/\log)$
of logarithmic differential
1-forms
is \'etale locally defined
by patching
$(\bigoplus_{i=0}^r
{\mathcal  O}_Xd\log T_i
\oplus
\bigoplus_{i=r+1}^n
{\mathcal  O}_XdT_i)/
(\sum_{i=0}^r
m_id\log T_i)$.
In the following,
we write
$$A^1_{X/S}
=
\Omega^1_{X/S}(\log/\log)
\quad\text{ and }
\quad
A^q_{X/S}=
\wedge^q_{{\mathcal  O}_X}
A^1_{X/S}$$ for short.
The ${\mathcal  O}_X$-module
$A^1_{X/S}$ is locally free
of rank $n=\dim X_K$.
The canonical map
$\Omega^1_{X/S}
\to
A^1_{X/S}$
induces an isomorphism
on the generic fiber $X_K$.
It also induces
an isomorphism
$\omega_{X/S}(
X_{s,{\rm red}}
-X_s)
\to
A^n_{X/S}$
where
$\omega_{X/S}
=\det\Omega^1_{X/S}$
denotes the
relative dualizing sheaf.

\begin{lm}\label{lmtm}
Let $X$ be
a proper generalized
semi-stable scheme
over $S={\rm Spec}\ {\mathcal  O}_K$.
Assume that
every irreducible component
of the closed fiber $X_s
=X\times_S{\rm Spec}\ F$
has odd multiplicity
in $X_s$.

{\rm 1.}
The inverse image
$I'\subset I$ of the pro-$2$ part
of ${\mathbb Z}_2(1)
\subset
\prod_{p'\neq p}
{\mathbb Z}_{p'}(1)$
acts trivially on
$H^q(X_{\bar K},
{\mathbb Q}_\ell)$
for every $q$.

{\rm 2.}
We define an effective
Cartier divisor $D$
by $2D=
X_s-
X_{s,{\rm red}}$.
Then,
the image of
$H^m(X,A^m_{X/S}(D))$
defines
a non-degenerate
${\mathcal  O}_K$-lattice
of 
$H^m(X_K,
\Omega^m_{X_K/K})$.
\end{lm}

\begin{proof}
1.
By the assumption,
the inverse image
$I'\subset I$ of the pro-$2$ part
acts trivially on
$R^q\psi
{\mathbb Q}_\ell$
for every $q$.
Hence, it follows
from the spectral sequence
$E_2^{p,q}=
H^p(X_{\bar F},
R^q\psi
{\mathbb Q}_\ell)\Rightarrow
H^{p+q}(X_{\bar K},
{\mathbb Q}_\ell)$.

2.
By the assumption
$X_s-X_{s, {\rm red}}=2D$,
we have an isomorphism
$\omega_{X/S}(-2D)
\to
A^n_{X/S}$.
It induces
an isomorphism
$A^m_{X/S}(D)
\to
{\mathcal  H}om
(A^m_{X/S}(D),
\omega_{X/S}).$
By Grothendieck duality,
it induces
an isomorphism
$$H^m(X,
A^m_{X/S}(D))/
(\text{torsion part})
\to
Hom_{{\mathcal  O}_K}
(H^m(X,A^m_{X/S}(D)),
{\mathcal  O}_K).$$
Hence the assertion follows.
\end{proof}

\begin{cor}\label{cortm}
Let $X$
be a proper
generalized semi-stable
scheme over 
$S={\rm Spec}\ {\mathcal  O}_K$
satisfying the condition in Lemma
{\rm \ref{lmtm}}.
Then, for a prime number
$\ell$ invertible in ${\mathcal  O}_K$,
we have
\setcounter{equation}0
\begin{equation}
\partial 
sw_2(H^n_\ell(X_K/K))=
\partial 
hw_2(H^n_{dR}(X_K/K))=0
\label{eqtm}
\end{equation}
in $H^1(F,{\mathbb Z}/2
{\mathbb Z})$.

In particular, further if the residue field
$F$ is finite,
Conjecture {\rm \ref{cn}}
is true.
\end{cor}

\begin{proof}
By Lemmas \ref{lmtm}.1
and \ref{lmp2}.1,
we have
$\partial 
sw_2(H^n_\ell(X_K/K))=0$.
By Lemmas \ref{lmtm}.2
and \ref{lmp2}.2,
we have
$\partial 
hw_2(H^m(X_K,\Omega^m_{X_K/K}))=0$.
Further by Lemma \ref{lmHdR},
we obtain
$\partial 
hw_2(H^n_{dR}(X_K/K))=0$.

Assume $F$ is finite.
Then, the map
$\partial\colon
H^2(K,{\mathbb Z}/2{\mathbb Z})
\to
H^1(F,{\mathbb Z}/2{\mathbb Z})$
is an isomorphism.
Since $2,d,-1$ are units in ${\mathcal  O}_K$,
the terms
$\{2,d\},\{d,-1\},\{-1,-1\}$ in 
(\ref{eqcn}) are 0.
Hence the equality (\ref{eqtm})
implies
Conjecture \ref{cn}.
\end{proof}

In the rest of this section,
we prove the assertion 1
in Theorem \ref{thm1}.

\begin{lm}\label{lmdvf}
Let $K$ be a complete discrete
valuation field
with residue field $F$
of characteristic
different from $2$.
Let $K'$
be a totally
ramified extension
of $K$ of degree $2$
and $I'
\subset G_{K'}$
be the inertia subgroup.
Let $\chi_{K'/K}
=d_{K'/K}
\in H^1(K,{\mathbb Z}/
2{\mathbb Z})
=K^\times/
K^{\times2}$ be the quadratic
character
of $G_K$ corresponding
to $K'$ or equivalently
the discriminant
of the quadratic extension $K'$.

{\rm 1.} 
Let $V$ be an orthogonal
representation of $G_K/I'
=G_F\times I/I'$
and let $V=V_0
\oplus (V_1\otimes
\chi_{K'/K})$ be
the decomposition
into the 
unramified part $V_0$
and the ramified
part $V_1\otimes
\chi_{K'/K}$.
Let $r$ be the degree of 
the representation
$V_1$
of
$G_F=G_K/I$.

Then we have
\setcounter{equation}0
\begin{equation}
\partial sw_2(V)=
\binom r2\{-1\}
+\det V_1
+\begin{cases}
0&\quad \text{ if $r$ is even}\\
\det V/\chi_{K'/K}
&\quad \text{ if $r$ is odd}.
\end{cases}
\label{eqdvf1}
\end{equation}

{\rm 2.} 
Let $(D,b)$ be a quadratic $K$-vector
space and $L$ be an ${\mathcal  O}_K$-lattice
of $D$ satisfying
${\mathfrak m}_KL^*\subset L\subset 
L^*=\{x\in D\mid b(x,y)\in {\mathcal  O}_K
\text{ for }y\in L\}$.
Then, there exists
a unique non-degenerate
${\mathcal  O}_{K'}$-lattice
$L'$ of $D_{K'}=
D\otimes_KK'$ satisfying
${\mathfrak m}_{K'}L'\subset 
{\mathcal  O}_{K'}L
\subset L'$.
The bilinear form $b$ 
induces a non-degenerate form on
the $F$-vector space
$\bar L_1=L'/{\mathcal O}_{K'}L$.
We put
$r=\dim \bar L_1$.

Then we have
\begin{equation}
\partial hw_2(D)=
\binom r2\{-1\}
+\disc \bar L_1
+\begin{cases}
0
&\quad \text{ if $r$ is even}\\
\disc D/d_{K'/K}
&\quad \text{ if $r$ is odd}.
\end{cases}
\label{eqdvf2}
\end{equation}
\end{lm}

\begin{proof}
1. 
Since $sw(V)=sw(V_0)\cdot 
sw(V_1\otimes
\chi_{K'/K})$,
we have
$$\partial sw_2(V)=
\partial sw_1(V_1\otimes
\chi_{K'/K})\cdot
sw_1(V_0)
+\partial sw_2(V_1\otimes
\chi_{K'/K}).$$
We have
$\partial sw_1(V_1\otimes
\chi_{K'/K})=r
\in H^0(F,{\mathbb Z}/
2{\mathbb Z})=
{\mathbb Z}/
2{\mathbb Z}$.
By Corollary \ref{cort2},
we have
$\partial 
sw_2(V_1\otimes 
\chi_{K'/K})
=(r-1)\det (V_1)
+\binom r2 \{-1\}$
since $\partial \chi_{K'/K}=1$
in $H^0(F,{\mathbb Z}/
2{\mathbb Z})=
{\mathbb Z}/
2{\mathbb Z}$
and $\chi_{K'/K}
\cup\chi_{K'/K}
=\chi_{K'/K}
\cup \{-1\}$
in $H^2(K,{\mathbb Z}/
2{\mathbb Z})$.
If $r$ is odd,
we have
$sw_1(V_0)
=\det V/\det V_1\cdot 
\chi_{K'/K}$.
Hence the assertion follows.

2. 
Since the canonical map
$\bar L_0=
L/{\mathfrak m}_KL^*
\to
L^*/{\mathfrak m}_KL^*
=Hom_F(L/{\mathfrak m}_KL,F)$
is injective,
$b$ induces a non-degenerate form on
$\bar L_0$.
Take a direct summand
$L_0$ of $L$ such that
$L_0/{\mathfrak m}_KL_0=
\bar L_0$
and put $L_1=
L_0^\perp
=\{x\in L\mid
b(x,y)=0\text{ for }
y\in L_0\}$.
Then, 
$L_0$ is non-degenerate and
we have
$L=L_0\oplus L_1$
and
$L^*=L_0\oplus {\mathfrak m}_K^{-1}L_1$.

The conditions
${\mathfrak m}_{K'}L'\subset 
{\mathcal  O}_{K'}L
\subset L'$
and that
$L'$ is non-degenerate
imply
$L'\subset 
{\mathcal  O}_{K'}L^*
\subset {\mathfrak m}_{K'}^{-1}L'$
and hence
${\mathcal  O}_{K'}L
+{\mathfrak m}_{K'}L^*
\subset L'
\subset
{\mathcal  O}_{K'}L^*
\cap{\mathfrak m}_{K'}^{-1}L$.
Since 
${\mathcal  O}_{K'}L
+{\mathfrak m}_{K'}L^*$
and ${\mathcal  O}_{K'}L^*
\cap {\mathfrak m}_{K'}^{-1}L$
are both equal to
${\mathcal  O}_{K'}L_0
\oplus
{\mathfrak m}_{K'}^{-1}L_1$,
the uniqueness of $L'$
follows.
It is clear that
$L'=
{\mathcal  O}_{K'}L_0
\oplus
{\mathfrak m}_{K'}^{-1}L_1$
is non-degenerate.
Since $L'/
{\mathfrak m}_{K'}L'$
is the orthogonal direct sum
$\bar L_0\oplus \bar L_1$,
the bilinear form $b$
induces non-degenerate forms on 
$\bar L_0$ and $\bar L_1$.

We put $D_0=KL_0$
and $D_1=KL_1$.
Since $hw(D)=hw(D_0)\cdot hw(D_1)$,
we have
$$\partial hw_2(D)=
\partial hw_1(D_1)\cdot
hw_1(D_0)
+\partial hw_2(D_1).$$
We have
$\partial hw_1(D_1)=r
\in H^0(F,{\mathbb Z}/
2{\mathbb Z})=
{\mathbb Z}/
2{\mathbb Z}$.
Since $L_1$
is non-degenerate with
respect to the
restriction of $\pi^{-1}b$
and since $\{\pi,\pi\}=
\{\pi,-1\}$,
we have
$\partial hw_2(D_1)
=(r-1)\partial \{\disc D_1,\pi\}
+\binom r2 \{-1\}$
by Lemma \ref{lmitp}.2.
If $r$ is even,
we have
$\partial \{\disc D_1,\pi\}
=
\disc \bar L_1$.
If $r$ is odd,
we have $hw_1(D_0)
=\disc D/\disc D_1$
and 
$\disc D_1=
\disc \bar L_1
\cdot d_{K'/K}$.
Thus the assertion follows
\end{proof}

The following Lemma
is inspired by
\cite{PL}.

\begin{lm}\label{lmDbl}
Let $X$ be a proper
generalized
semi-stable scheme
over 
$S={\rm Spec}\ {\mathcal  O}_K$.
Assume
that the divisor
$D=X_s-
X_{s,{\rm red}}$
is smooth over
$F$
and that
the complement
$X\setminus D$
is smooth over $S$.
Let $K'$
be a totally ramified
extension of
degree $2$ and
put
$S'={\rm Spec}\
{\mathcal  O}_{K'}$.

{\rm 1.}
The normalization $X'$
of the base change
$X\times_SS'$ is
semi-stable.
The reduced inverse image
$D'=(D\times_XX')_{\rm red}$
is smooth over $F$.
The double covering
$\varphi\colon
D'\to D$ is
\'etale outside
$C=D\cap (X_s-2D)$
and is totally ramified
along $C$.

{\rm 2.}
Let $\ell$ be a
prime number invertible
in ${\mathcal  O}_K$ and
$I'\subset G_{K'}\subset G_K$
be the inertia subgroup.
Then, the Galois representation
on $H^q(X_{\bar K},
{\mathbb Q}_\ell)$
factors through
the quotient
$G_K/I'$.

Let 
$H^q(X_{\bar K},
{\mathbb Q}_\ell)^-
\subset
H^q(X_{\bar K},
{\mathbb Q}_\ell)$
and
$H^q(D'_{\bar F},
{\mathbb Q}_\ell)^-
\subset
H^q(D'_{\bar F},
{\mathbb Q}_\ell)$
be the minus-parts
with respect to
the actions of
$I/I'$.
Then, the cospecialization map
$H^q(X'_{\bar F},
{\mathbb Q}_\ell)
\to
H^q(X_{\bar K},
{\mathbb Q}_\ell)$
induces
an isomorphism
$H^q(D'_{\bar F},
{\mathbb Q}_\ell)^-
\to
H^q(X_{\bar K},
{\mathbb Q}_\ell)^-$.

{\rm 3.}
Assume that $X_K$
is of
even dimension $n=2m$.
Then, the canonical map

\noindent
$H^m(X',A^m_{X'/S'})
\otimes_{{\mathcal  O}_{K'}}F
\to 
H^m(X'_F,A^m_{X'_F/F})$
is injective
and the image
is the orthogonal
of
the image
of the torsion part
$H^m(X',A^m_{X'/S'})_{\rm tors}$.

Let
$H^m(X'_F,A^m_{X'_F/F})=
H^m(X'_F,A^m_{X'_F/F})^+
\oplus
H^m(X'_F,A^m_{X'_F/F})^-$
be the decomposition
in to the plus part
and the minus part
with respect
to the action of
$I/I'={\rm Gal}(K'/K)$.
Then, the image
of
$H^m(X,A^m_{X/S})
\otimes_{{\mathcal  O}_K}F$
in the plus-part
$H^m(X'_F,A^m_{X'_F/F})^+$
is equal to
the image
of
$H^m(X',A^m_{X'/S'})
\otimes_{{\mathcal  O}_{K'}}F$
and the image
of
$H^m(X,A^m_{X/S})
\otimes_{{\mathcal  O}_K}F$
in the minus-part
$H^m(X'_F,A^m_{X'_F/F})^-
=
H^m(D',\Omega^m_{D'/F})^-$
is $0$.
\end{lm}

\begin{proof}
1.
Since the assertion
is \'etale local
on $X$,
we may assume
$X$ is smooth over
$A={\mathcal  O}_K[x,y]/(xy^2-\pi)$
and $K'=K(\sqrt{\pi})$
for a prime element
$\pi$ of $K$.
Then, the normalization
of
$A\otimes_{
{\mathcal  O}_K}{\mathcal  O}_{K'}$
is
${\mathcal  O}_{K'}[x',y]/
(x'y-\sqrt{\pi})$
where $x=x^{\prime 2}$
and the assertion follows.

2.
We consider the
spectral sequence
$E_2^{p,q}=
H^p(X_{\bar F},
R^q\psi
{\mathbb Q}_\ell)\Rightarrow
H^{p+q}(X_{\bar K},
{\mathbb Q}_\ell)$.
We recall
computation
of $R^q\psi
{\mathbb Q}_\ell$
from \cite[Proposition 6]{St}.
We have
$R^q\psi{\mathbb Z}_\ell=0$
for $q>1$.
Further,
the restriction
$R^0\psi{\mathbb Z}_\ell|_D$
is canonically
isomorphic
to $\varphi_*{\mathbb Z}_{\ell,D'}$
and
the restriction
$R^0\psi
{\mathbb Z}_\ell|_{X_s-2D}$
is canonically
isomorphic to
${\mathbb Z}_\ell|_{X_s-2D}$.
The isomorphism
is compatible
with the actions of $I$.
The sheaf
$R^1\psi{\mathbb Z}_\ell$
is the direct
image of
a locally constant
sheaf of rank 1
on the intersection
$C=D\cap (X_s-2D)$
with the trivial
action of $I$.

Hence,
the inertia action on
$H^p(X_{\bar F},
R^q\psi
{\mathbb Q}_\ell)$
is trivial for $q\neq 0$.
For $q=0$,
the action of
$I'$ on
$H^p(X_{\bar F},
R^0\psi
{\mathbb Q}_\ell)$
is trivial.
Further, the minus part
of 
$H^p(X_{\bar F},
R^0\psi
{\mathbb Q}_\ell)$
is isomorphic to 
$H^p(D'_{\bar F},
{\mathbb Q}_\ell)^-.$

3.
We consider the 
commutative diagram
$$\begin{CD}
0
\to&
H^m(X',A^m_{X'/S'})
\otimes_{{\mathcal  O}_{K'}} F
@>>>
H^m(X'_F,A^m_{X'_F/F})
@>>>
Tor^{{\mathcal  O}_{K'}}_1
(H^{m+1}(X',A^m_{X'/S'}),
F)&
\to 0\\
&@.@VVV @VVV\\
&@.
H^m(X'_F,A^m_{X'_F/F})^\vee
@>>>
(
H^m(X',A^m_{X'/S'})_{\rm tors}
\otimes_{{\mathcal  O}_{K'}}F)^\vee.
\end{CD}
$$
The upper line
is exact
and
$^\vee$
denote
the $F$-linear dual.
By Grothendieck duality,
the vertical arrows
are isomorphisms.
Thus, the first
paragraph is proved.

We consider
the decomposition
${\mathcal  O}_{D'}
=
{\mathcal  O}_D\oplus
{\mathcal  L}$
by the action
of the Galois group
${\rm Gal}(K'/K)
={\rm Gal}(D'/D)$.
The computation 
in the proof of 1.\
shows
that we also have a
decomposition
${\mathcal  O}_{X'_F}
=
{\mathcal  O}_{X_{F,{\rm red}}}
\oplus
{\mathcal  L}$.
We put $A^m_{D/F}=
\Omega^m_{D/F}(\log C)$
and $A^m_{D'/F}=
\Omega^m_{D'/F}(\log C)$.
The canonical map
$A^m_{D/F}\otimes_
{{\mathcal  O}_D}
{\mathcal  O}_{D'}\to
A^m_{D'/F}
=
A^m_{D/F}\oplus
(A^m_{D/F}\otimes_
{{\mathcal  O}_D}{\mathcal  L})$
is an isomorphism.
Since the sequences
$0\to
\Omega^m_{D/F}\to
A^m_{D/F}\to
\Omega^{m-1}_{C/F}\to
0$
and
$0\to
\Omega^m_{D'/F}\to
A^m_{D'/F}\to
\Omega^{m-1}_{C/F}\to
0$ are exact,
the inclusion
$\Omega^m_{D'/F}
\to
A^m_{D'/F}$
induces an isomorphism
$\Omega^{m-}_{D'/F}
\to
A^{m-}_{D'/F}=
A^m_{D/F}\otimes_
{{\mathcal  O}_D} {\mathcal  L}$
on the minus parts.
Since the canonical map
$A^m_{X'/S'}
\otimes _{{\mathcal  O}_{X'}}
{\mathcal  O}_{D'}
\to
A^m_{D'/F}$
is an isomorphism,
it is also identified with
$A^{m-}_{X'_F/F}$.
Thus, we have
$H^m(X'_F,A^m_{X'_F/F})^-
=
H^m(D',\Omega^m_{D'/F})^-$.

Further
the computation 
in the proof of 1.\
shows that
the canonical map
${\rm Coker}(
{\mathcal  O}_{X\times_SS'}
\to
{\mathcal  O}_{X'})
\to
{\rm Coker}(
{\mathcal  O}_D
\to
{\mathcal  O}_{D'})={\mathcal  L}$
is an isomorphism.
Hence, we have exact
sequences
$0\to A^m_{X/S}
\otimes_{{\mathcal  O}_X}
{\mathcal  O}_{X\times_SS'}
\to
A^m_{X'/S'}
\to
\Omega^{m-}_{D'/F}
\to 0$
and
$$\to
H^m(X,A^m_{X/S})
\otimes_{{\mathcal  O}_K}
{\mathcal  O}_{K'}
\to
H^m(X',A^m_{X'/S'})
\to
H^m(D',\Omega^m_{D'/F})^-
\to.
$$
Thus, we obtain
an exact sequence
$$\to
H^m(X,A^m_{X/S})
\otimes_{{\mathcal  O}_K}F
\to
H^m(X',A^m_{X'/S'})
\otimes_{{\mathcal  O}_{K'}}F
\to
H^m(D',\Omega^m_{D'/F})^-
\to
$$
and the assertion
in the second paragraph follows.
\end{proof}

\begin{pr}\label{prdbl}
Let the notation be as in
Lemma {\rm \ref{lmDbl}}.

{\rm 1.}
The following conditions
are equivalent:

{\rm (1)}
The determinant
$\det H^n(X_{\bar K},
{\mathbb Q}_\ell)$ is unramified.

{\rm (2)}
The discriminant
$d_X=
\disc H^n_{dR}(X/K)
\in K^\times/
K^{\times2}$
is in the image of
$F^\times/
F^{\times2}.$

{\rm (3)}
The Euler number
$\chi(D_{\bar F},
{\mathbb Q}_\ell)
=\chi_{dR}(D/F)$
is even.

{\rm 2.}
Let $\chi^-(D')$
denote the Euler number
$\sum_q(-1)^q
\dim H^q(D'_{\bar F},
{\mathbb Q}_\ell)^-$
of the minus-part
and let
$\chi_{K'/K}
\in H^1(K,{\mathbb Z}/
2{\mathbb Z})$ denote
the character of order
$2$ corresponding
to the quadratic extension 
$K'$ over $K$.
Then,
the image of the left hand
side $sw_2
(H^n_\ell(X))
+\{e,-1\}+\beta\cdot c_\ell$
of {\rm (\ref{eqcn})}
by $\partial$
is equal to
\setcounter{equation}0
\begin{align}
&\binom
{\chi^-(D')}2
\{-1\}+
\det \left(H^n(D'_{\bar F},
{\mathbb Q}_\ell)^-
\left(\frac n2\right)\right)
\label{eqV-}
\\
&\qquad
+\begin{cases}
0
&\quad
\text{ if 
$\det H^n(X_{\bar K},
{\mathbb Q}_\ell)$ is unramified},\\
\det H^n(X_{\bar K},
{\mathbb Q}_\ell)
/\chi_{K'/K}
&\quad
\text{ if 
$\det H^n(X_{\bar K},
{\mathbb Q}_\ell)$ is ramified}
\end{cases}
\nonumber
\end{align}
in $H^1(F,{\mathbb Z}/
2{\mathbb Z})$.

{\rm 3.}
Let $h^{m,m-}(D')$
denote the dimension
$\dim H^m(D'_{\bar F},
\Omega^m_{D'/F})^-$
of the minus-part
and let
$d_{K'/K}
\in H^1(K,{\mathbb Z}/
2{\mathbb Z})$ denote
the discriminant
of the quadratic extension 
$K'$ over $K$.
Let $H^n_{dR}(D'/F)'$
denote $H^n_{dR}(D'/F)$
with the symmetric bilinear
form multiplied by $(-1)^{n/2}$
and let
$H^n_{dR}(D'/F)^{\prime-}$
denote its minus part
with respect to the
action of ${\rm Gal}(D'/D)$.
Let $r$ be as in {\rm(\ref{eqr})}.

Then,
the image of the right hand
side of {\rm (\ref{eqcn})}
by $\partial$
is equal to
\begin{align}
\binom
{b_{n,dR}^-(D')}2
\{-1\}+
&\disc H^n_{dR}(D'/F)^{\prime-}
+\chi_{dR}(D/F)\cdot\{2\}
\label{eqD-}
\\
&\quad+
\begin{cases}
0&\quad
\text{ if }
d_X\in F^\times/
F^{\times2},\\
(d'_X-d_{K'/K})
+
r\{-1\}
&\quad
\text{ if }
d_X
\notin F^\times/
F^{\times2}
\end{cases}
\nonumber
\end{align}
in $H^1(F,{\mathbb Z}/
2{\mathbb Z})$.
\end{pr}

\begin{proof}
1.
By Lemma \ref{lmDbl}.2,
the condition (1)
is equivalent to
that
$\dim H^n(D'_{\bar F},
{\mathbb Q}_\ell)^-$
is even.
Hence, it
is equivalent to
that the Euler number
$\chi(D'_{\bar F},
{\mathbb Q}_\ell)^-$
is even.
Since $\chi(D'_{\bar F},
{\mathbb Q}_\ell)^-
=\chi(D_{\bar F},
{\mathbb Q}_\ell)
-
\chi(C_{\bar F},
{\mathbb Q}_\ell)$
and since
the Euler number
$\chi(C_{\bar F},
{\mathbb Q}_\ell)$
is even for
proper smooth
scheme $C_{\bar F}$
of odd dimension,
it is equivalent to
the condition (3).

By Lemma \ref{lmDbl}.3,
the condition (2)
is equivalent to
that
$\dim H^m(D',
\Omega^m_{D'/F})^-$
is even.
By Serre duality
and the Hodge to de Rham
spectral sequence,
it is equivalent to
that the Euler number
$\chi_{dR}(D'/F)^-$
is even.
Since $\chi_{dR}(D'/F)^-
=\chi_{dR}(D/F)
-
\chi_{dR}(C/F)$
and since
the Euler number
$\chi_{dR}(C/F)$
is even for
proper smooth
scheme $C$
of odd dimension,
it is also equivalent to
the condition (3).

2.
We put 
$b_n^-=
\dim H^n(D'_{\bar F},
{\mathbb Q}_\ell)^-$.
By Lemmas \ref{lmdvf}.1
and \ref{lmDbl}.2,
we have
\begin{align*}
\partial(sw_2(H^n_\ell(X)))
=&\
\binom
{b_n^-}2\{-1\}
+\det \left(
H^n(D'_{\bar F},
{\mathbb Q}_\ell)^-
\left(\frac n2\right)\right)\\
&\quad
+\begin{cases}
0
&\quad \text{ if $b_n^-$ is even,}\\
\det H^n(X_{\bar K},
{\mathbb Q}_\ell)^-\left(\frac n2\right)
/d_{K'/K}
&\quad \text{ if $b_n^-$ is odd}.
\end{cases}
\end{align*}
Further by Lemma \ref{lmDbl}.2,
the character
$\det H^n(X_{\bar K},
{\mathbb Q}_\ell)^-$
is unramified
if and only if
$b_n^-$ is even.
We put
$r^-=\sum_{q<n}
(-1)^q\dim H^q(D'_{\bar F},
{\mathbb Q}_\ell)^-$.
Then, 
further
by Lemma \ref{lmDbl}.2,
the character
$e=
\prod_{q<n}
(\det H^q(X_{\bar K},
{\mathbb Q}_\ell))^{(-1)^q}$
is unramified
if and only
if $r^-$ is even.
Hence, we have
$\partial \{e,-1\}
=r^-\{-1\}$.
Since $\chi^-(D')
=b_n^-+2r^-$,
we have
$\displaystyle{
\binom{\chi^-(D')}2
=\binom{b_n^-}2+r^-}$.
Since 
$\partial(c_\ell)=0$,
the assertion follows.

3.
By (\ref{eqcnhg}),
the right hand side of
(\ref{eqcn})
is equal to
$hw_2(H^\bullet_{dR}
(X_K/K))
+\{2,d_X\}
+\eta(c_\ell-c_2)$
as in Corollary \ref{corg}.
Let $H^m(X_K,\Omega^m_{X/K})'$
denote the 
symmetric $K$-vector space
$H^m(X_K,\Omega^m_{X/K})$
with the symmetric bilinear
form multiplied by $(-1)^m$.
We put
$h^{m,m}=
H^m(X_K,\Omega^m_{X/K})$
and define $r'$ by
$\chi_{DR}(X/K)
=h^{m,m}+2r'$.
Then, by Lemma \ref{lmHdR},
we have
\begin{align*}
&\partial(
hw_2(H^\bullet_{dR}
(X_K/K)))\\
=&\
\partial(
hw_2(H^m(X_K,\Omega^m_{X/K})'))
+\begin{cases}
0&\quad \text{ if }
\disc H^m(X_K,\Omega^m_{X/K})
\in F^\times/
F^{\times2}
\\
r'\cdot
\{-1\}
&\quad \text{ if }
\disc H^m(X_K,\Omega^m_{X/K})
\notin F^\times/
F^{\times2}.
\end{cases}
\end{align*}
The condition
$\disc H^m(X_K,\Omega^m_{X/K})
\in F^\times/
F^{\times2}$
is equivalent to
$d_X=\disc H^n_{dR}(X/K)
\in F^\times/
F^{\times2}$.
It is further equivalent to
that 
$h^{m,m-}(D')$ is even.

We apply
Lemma \ref{lmdvf}.2.\
to $H^m(X_K,\Omega^m_{X/K})'$
and define $\bar L_1$
as in Lemma \ref{lmdvf}.2.
Then, we obtain
\begin{align*}
\partial(
hw_2(H^m(X_K,\Omega^m_{X/K})'))
=&\
\binom
{\dim \bar L_1}2\{-1\}
+\disc \bar L_1\\
&+
\begin{cases}
0&\quad \text{ if $h^{m,m-}(D')$ is even}\\
\disc H^m(X_K,\Omega^m_{X_K/K})'
/d_{K'/K}
&\quad \text{ if $h^{m,m-}(D')$ is odd}.
\end{cases}
\end{align*}
We put
$s=\dim
H^m(X',A^m_{X'/S'})_{\rm tors}
\otimes_{{\mathcal  O}_{K'}}F$.
Then, we have
$h^{m,m-}(D')=
\dim \bar L_1 +
2\cdot s$
and
$\disc H^m(D',
\Omega_{D'/F})^{\prime -}
=
\disc \bar L_1
+s\cdot \{-1\}$
by Lemma \ref{lmDbl}.3.
Hence, for the right hand side,
we have
$$
\binom
{\dim \bar L_1}2\{-1\}
+\disc \bar L_1
=
\binom
{h^{m,m-}(D')}2\{-1\}
+\disc H^m(D',
\Omega_{D'/F})^{\prime -}.$$
It is further 
equal to
$\binom
{b_{n,dR}(D')^-}2\{-1\}
+\disc H^n_{dR}
(D'/F)^{\prime -}$
by Lemma \ref{lmHdR}.

We have
$\partial\{d,2\}
=
h^{m,m-}(D')
\cdot\{2\}$
and
$h^{m,m-}(D')
\equiv \chi^-(D')\bmod 2$.
Further, as
in the end of the proof of 1.,
we have
$\chi^-(D')
\equiv \chi_{dR}(D/F)\bmod 2$.
We also have
$\partial(c_\ell)=
\partial(c_2)=0$.
Since 
$\disc H^m(X_K,\Omega^m_{X/K})'
-
\disc H^n_{dR}(X_K/K)'
=
(r-r')\{-1\}$,
the assertion follows.
\end{proof}

\begin{lm}\label{lmdbl}
Let $D$ be
a projective smooth
scheme of even dimension $n$
over a field $F$
of characteristic $\neq 2$.
Let $D'\to D$
be a double covering
ramified
along a smooth
divisor $C\subset D$.
Let
$H^n_{dR}(D'/F)'$
be the quadratic $F$-vector space
$H^n_{dR}(D'/F)$
with symmetric bilinear
form multiplied by $(-1)^{n/2}$
and let 
$H^n_{dR}(D'/F)^{\prime -}$
denote the minus part
with respect to the action of
${\rm Gal}(D'/D)=\{\pm1\}$.
We put $r_{D'}^-
=\sum_{q<n}(-1)^q
\dim 
H^q_{dR}(D'/F)^-$.

Then, we have
\setcounter{equation}0
\begin{equation}
\det
H^n(D'_{\bar F},
{\mathbb Q}_\ell)^-
=
\disc 
H^n_{dR}(D'/F)^{\prime -}
+r_{D'}^-\cdot\{-1\}+
\chi_{dR}(D/F)
\cdot \{2\}
\end{equation}
in $H^1(F,{\mathbb Z}/
2{\mathbb Z})$.
\end{lm}

\begin{proof}
We put $r_D
=\sum_{q<n}(-1)^q
\dim 
H^q_{dR}(D/F)$ and
$r_{D'}
=\sum_{q<n}(-1)^q
\dim 
H^q_{dR}(D'/F)$.
Let $H^n_{dR}(D/F)'$
denote the quadratic $F$-vector space
$H^n_{dR}(D/F)$
with symmetric bilinear
form multiplied by $(-1)^{n/2}$.
Then, by \cite[Theorem 2]{S1},
we have
$\det
H^n(D_{\bar F},
{\mathbb Q}_\ell)
=
\disc 
H^n_{dR}(D/F)'
+r_D\{-1\}$
and
$\det
H^n(D'_{\bar F},
{\mathbb Q}_\ell)
=
\disc 
H^n_{dR}(D'/F)'
+r_{D'}\{-1\}$.
On the left hand sides,
we have
$$\det H^n(D'_{\bar F},
{\mathbb Q}_\ell)
=
\det H^n(D_{\bar F},
{\mathbb Q}_\ell)
\cdot
\det H^n(D'_{\bar F},
{\mathbb Q}_\ell)^-$$
since 
$H^n(D'_{\bar F},
{\mathbb Q}_\ell)
=
H^n(D_{\bar F},
{\mathbb Q}_\ell)
\oplus
H^n(D'_{\bar F},
{\mathbb Q}_\ell)^-$.
On the other hand,
since 
$H^q_{dR}(D'/F)
=
H^q_{dR}(D/F)
\oplus
H^q_{dR}(D'/F)^-$,
we have
$r_{D'}=r_D+r_{D'}^-$.
Further since
the restriction of
the symmetric bilinear
form of
$H^n_{dR}(D'/F)$
on
$H^n_{dR}(D/F)
\subset
H^n_{dR}(D'/F)$
is $2$-times
that of
$H^n_{dR}(D/F)$,
we have
$$\disc H^n_{dR}(D'/F)'
=
2^{\chi(D)}
\cdot
\disc H^n_{dR}(D/F)'
\cdot
\disc H^n_{dR}(D'/F)^{\prime-}.$$
Hence the assertion follows.
\end{proof}

\begin{cor}
\label{cordvf}
Let the notation
be as in Lemma {\rm \ref{lmDbl}}.
Assume $D$ and $X_K$
are projective.
Then, the both sides
of {\rm (\ref{eqcn})}
have the same
images by $\partial$.
In particular,
further if the residue
field $F$ is finite,
Conjecture {\rm \ref{cn}}
is true in this case.
\end{cor}

\begin{proof}
We compare
the sums of the terms
in the first lines
in (\ref{eqV-})
and (\ref{eqD-}).
Since
$\chi^-(D')=
b_{n,dR}^-(D')
+2r_{D'}^-$,
they are equal to each other
by Lemma \ref{lmdbl}.
In the case
where the equivalent conditions
in Proposition \ref{prdbl}.1
do not hold,
the remaining terms
are also equal to each other
by \cite[Theorem 2]{S1}.
\end{proof}

\begin{cor}\label{corodp}
Let $K$ be a complete
discrete valuation
field with
residue field
$F$ of characteristic
$\neq 2,\ell$.
Let $X_{{\mathcal  O}_K}$ be
a proper regular flat scheme
over a discrete valuation ring ${\mathcal  O}_K$.
Assume that the generic fiber
$X_K$ is projective
and smooth of even dimension
and that the closed fiber
$X_s$ has
at most ordinary double points
as singularities.

Then, the both sides
of {\rm (\ref{eqcn})}
have the same
images by $\partial$.
In particular,
further if the residue
field $F$ is finite,
Conjecture {\rm \ref{cn}}
is true in this case.
\end{cor}

\begin{proof}
The blow-up of $X_{{\mathcal  O}_K}$
at the singular points
of the closed fibers
satisfies the assumption
of Corollary \ref{cordvf}.
\end{proof}

\section{Cristalline representations}\label{scris}

In this short section,
we derive the assertion 2
in Theorem \ref{thm1}
from the following computation
of the second Stiefel-Whitney
class of an orthogonal
cristalline representation.

\begin{lm}[{{\rm \cite[Theorem 2.3]{S2}}}]\label{lm3}
Let $K$ be a complete
discrete valuation field
of characteristic $0$
with perfect residue field
$F$ of characteristic $p>2$.
We assume that
$p$ is a prime element of $K$.
Let $V$ be an orthogonal
cristalline ${\mathbb Q}_p$-representation
of the absolute Galois group 
$G_K
=\text{\rm Gal}(\bar K/K)$
and $D=D_{\rm cris}(V)$
be the associated
quadratic $K$-vector space.
Assume that
the Hodge filtration
$F^\bullet$ on $D$ satisfies
$F^{\frac {p-1}2}D=0$.

Then we have 
$$\partial
sw_2(V)
= \sum_{q>0}
q\cdot
\dim_KGr_F^qD\cdot
\partial c_p$$
in $H^1(F,{\mathbb Z}/2)$.
In particular,
further if the residue field
$F$ is finite,
we have
$sw_2(V)
= \sum_{q>0}
q\cdot
\dim_KGr_F^qD\cdot
c_p$
in $H^2(K,{\mathbb Z}/2)$.
\end{lm}

Applying Lemma \ref{lm3}, we 
compute the second Stiefel-Whitney class
in a good reduction case.

\begin{pr}\label{pr3}
Let $K$ be a complete
discrete valuation field
of characteristic $0$
with perfect residue field
$F$ of characteristic $p\neq 0$.
We assume that
$p$ is a prime element of $K$.
Let $X$ be a proper smooth scheme
of even dimension $n$
over $K$
and assume that
$X$ has a proper smooth model 
$X_{{\mathcal  O}_K}$ over
the integer ring ${\mathcal  O}_K$.
We further assume that
$H^q(X,\Omega^{n-q}_{X/K})=0$
for $|q- \frac n2|\ge \frac {p-1}2$.

We put $h=\sum_{q<\frac n2}
\left(\frac n2-q\right)
\dim H^{n-q}(X,\Omega^q_{X/K})$ 
{\rm (\ref{eqh})} as in 
Lemma {\rm \ref{lmch0}}.
Then we have
$$\partial 
sw_2(H^n_p(X))=h\cdot 
\partial c_p$$
in $H^1(F,{\mathbb Z}/2)$.
In particular,
further if $K$ is a finite
extension of ${\mathbb Q}_p$,
we have
$sw_2(H^n_p(X))=h\cdot 
 c_p$
in $H^2(K,{\mathbb Z}/2)$.
\end{pr}

\begin{proof}
Since $H^n(X_{\bar K},{\mathbb Q}_p)$
is a cristalline representation
by \cite{FM},
it is enough to apply Lemma \ref{lm3}
if $p>2$.
If $p=2$,
the assumption
implies that
$H^n(X_{\bar K},{\mathbb Q}_p)(\frac n2)$
is unramified
and $h=0$.
Hence the assertion follows.
\end{proof}

\begin{cor}\label{cor3}
We keep the assumption 
in Proposition {\rm \ref{pr3}}.
The both sides
of {\rm (\ref{eqcn})}
have the same
images by $\partial$.
In particular,
further if $K$ is a finite
unramified
extension of ${\mathbb Q}_p$,
Conjecture {\rm \ref{cn}}
is true in this case.
\end{cor}

The assumption 
$H^q(X,\Omega^{n-q}_{X/K})=0$
for $|q- \frac n2|\ge \frac {p-1}2$
of Proposition \ref{pr3}
is satisfied if
$n+1< p=\ell$.
Hence, 
Corollary \ref{cor3}
implies the assertion 2
in Theorem \ref{thm1}
for $p=\ell\neq2$.
The case $p=\ell=2$
will be proved
as a consequence
of Proposition \ref{prhp}.

\begin{proof}
By Corollary \ref{cortm},
we have
$\partial hw_2(H^n_{dR}(X/K))=0$.
By Lemma \ref{lmtm},
we also have
$\partial d_X=0$
and 
$\partial e=0$.
Thus, the image by $\partial$
of the left hand side
of (\ref{eqp}) is
$\partial sw_2(H^n_\ell(X))$
and that of the right hand
side
is $h\cdot \partial c_\ell$.
Hence, it follows
from Proposition \ref{pr3}.
\end{proof}

\section{Hodge structures}\label{sHdg}

In this section, we prove 
the assertion 3 of
Theorem \ref{thm1}
for a projective smooth variety over
${\mathbb R}$,
using polarizations of Hodge structures.

Before starting proof, we
recall some terminology 
on Hodge structures.
An ${\mathbb R}$-Hodge structure of weight 0
over ${\mathbb R}$ is
an ${\mathbb R}$-vector space $V$
of finite dimension endowed with
the following structures:
A representation $\text{Gal}({\mathbb C}/{\mathbb R})
\to GL_{\mathbb R}(V)$ and
a decreasing filtration $F^\bullet$ on
the ${\mathbb R}$-vector space
$D=(V\otimes_{\mathbb R}{\mathbb C})^{\text{Gal}({\mathbb C}/{\mathbb R})}$
giving a Hodge decomposition
$V_{\mathbb C}=V\otimes_{\mathbb R}{\mathbb C}=
\bigoplus_pV^{p,-p}$ where $V^{p,q}=F_{\mathbb C}^p\cap
\overline{F_{\mathbb C}^q}$.
Here $\sigma\in \text{Gal}({\mathbb C}/{\mathbb R})$ acts on
$V_{\mathbb C}$ as $\sigma\otimes \sigma$
and ${ }^-$ denotes 
$\sigma\otimes 1$.
We say an ${\mathbb R}$-Hodge structure of weight 0
over ${\mathbb R}$ is polarized if it is equipped with
a non-degenerate symmetric bilinear form
$b_V:V\times V\to{\mathbb R}$ satisfying the following condition.
Let $C$ be an automorphism of $V_{\mathbb C}$ which is
defined to be
the multiplication by $i^{-p}\bar i^{-q}=(-1)^p$ on
$V^{p,q}=V^{p,-p}$. Then the condition is that
the bilinear form $(x,y)\mapsto b_V(x,Cy)$ on $V$
is symmetric and positive definite.
Let $b_D$ denote the symmetric bilinear form
on $D$ defined as 
the restriction of
the symmetric bilinear form on $V_{\mathbb C}$
induced by $b_V$.
We put $h^{p,q}=\dim_{\mathbb C}V^{p,q}$
and let $h^{0,0,\pm}$ be the
multiplicity of eigenvalue $\pm1$
of the complex conjugation
on $V^{0,0}_{\mathbb R}=V\cap V^{0,0}$.

For an orthogonal representation $V$
of $\text{Gal}({\mathbb C}/{\mathbb R})$,
let $v^-=\dim V^-$
be the dimension of the $(-1)$-eigenspace
$V^-=\{x\in V\mid 
\sigma(x)=-x\}$
where $\sigma\in \text{Gal}({\mathbb C}/{\mathbb R})$
denotes the complex conjugate.
Then we have
$sw_1(V)=v^-\{-1\}$ and
$sw_2(V)=\binom{v^-}2\{-1,-1\}$.
For a quadratic ${\mathbb R}$-vector space $(D,b_D)$,
let $d^-=\dim D^-$ be the dimension of
a maximal subspace $D^-$ of
$D$ where the symmetric bilinear form
$b_D|_{D^-}$ is negative definite.
Then we have
$hw_1(D)=d^-\{-1\}$ and
$hw_2(D)=\binom{d^-}2\{-1,-1\}$.

\begin{lm}\label{lm22}
Let $V$ be a polarized ${\mathbb R}$-Hodge
structure of weight $0$ over ${\mathbb R}$.
Then we have
$$v^-=d^-=\sum_{p>0}h^{p,-p}+h^{0,0,-}.$$
\end{lm}

\begin{proof}
Since $V=
V^{0,0}_{\mathbb R}
\oplus 
(V\cap \bigoplus_{p>0}V^{p,-p})$,
it suffices to consider the cases
where
$V=V^{0,0}_{\mathbb R}$
and $V^{0,0}_{\mathbb R}=0$
respectively.
If $V=V^{0,0}_{\mathbb R}$,
the symmetric bilinear form
$b_V$ is positive definite and
we obtain
$v^-=d^-=h^{0,0,-}$ and 
$h^{p,-p}=0$ for $p>0.$
If $V^{0,0}_{\mathbb R}=0$,
we have $h^{0,0,-}=0$ and 
$v^-,\ d^-$ and $\sum_{p>0}h^{p,-p}$
are equal to $\dim V/2$.
\end{proof}


By the comparison theorem of
singular cohomology and
\'etale cohomology,
we have
$sw_i(H^n_\ell(X/{\mathbb R}))
=sw_i(H^n(X({\mathbb C}),{\mathbb Q}))$.
Now we are ready to
prove the following main
result of this section.

\begin{pr}\label{prHdg}
Let $X$ be a projective smooth scheme
over ${\mathbb R}$ of even dimension $n$.
Let 
$\beta$ {\rm (\ref{eqbeta})},
$r$ {\rm (\ref{eqr})} and 
$e=\sum_{q<n}e_q
\in H^1({\mathbb R},
{\mathbb Z}/
2{\mathbb Z})$ be 
as in Conjecture {\rm \ref{cn}}
and
let
$hw_2(H^n_{dR}(X/{\mathbb R})')
\in H^2({\mathbb R},
{\mathbb Z}/
2{\mathbb Z})$ and
$d'_X\in H^1({\mathbb R},
{\mathbb Z}/
2{\mathbb Z})$ 
be as in 
Lemma {\rm \ref{lm'}}.
Then we have
\setcounter{equation}0
\begin{align}
&\
sw_2(H^n_\ell(X/{\mathbb R}))
+\{e,-1\}
+\beta\cdot \{-1,-1\}
\label{eqHdg}\\
&\
=
hw_2(H^n_{dR}(X/{\mathbb R})')
+r\{d'_X,-1\}+
\binom r2
\{-1,-1\}
\nonumber
\end{align}
in $H^2({\mathbb R},{\mathbb Z}/2{\mathbb Z})$.
\end{pr}

Since 
$\{2,d\},c_\ell-c_2
\in
H^2({\mathbb R},{\mathbb Z}/2{\mathbb Z})$
in the formula (\ref{eqcn'})
are $0$,
Proposition \ref{prHdg}
implies the assertion 3
in Theorem \ref{thm1}.

\begin{proof}
We take an ample invertible sheaf and
consider
the associated Lefschetz decomposition
$H^n(X({\mathbb C}),{\mathbb R})(\frac n2)=
\bigoplus_{q\le n, q:\text{even}}
P^q$ by the primitive parts.
Then each $P^q$ is an ${\mathbb R}$-Hodge structure of
weight 0 over ${\mathbb R}$.
Further $(-1)^{\frac q2}$-times the restriction of $b_{\mathcal L}$
to $V^q$ defines a polarization on it.
Let $D=H^n_{dR}(X/{\mathbb R})=
\bigoplus_{q\le n, q:\text{even}}
D^q$ be the corresponding Lefschetz decomposition.

For an even integer
$0\le q\le n$,
let $e_q^-$ be
the multiplicity
of the eigenvalue
$-1$ of the action
of the complex conjugate
on $H^q(X_{\mathbb C},
{\mathbb Q})(\frac q2)$
and we put 
$e^-=
\sum_{q<n,\text{\rm even}}e_q^-$.
Let $(d^{\prime +},d^{\prime -})$
be the signature
of the quadratic vector
space $H^n_{dR}(X/{\mathbb R})'$
with the symmetric bilinear
form multiplied by $(-1)^{n/2}$.
We put
$P^+=
\bigoplus_{q< n,q\equiv 0\bmod 4}P^q$
and
$P^-=
\bigoplus_{q< n,q\equiv 2\bmod 4}P^q$
as in Lemma \ref{lmch0}.2 and
prove the congruence
\begin{equation}
e_n^-
-2e^-
\equiv
d^{\prime -}
-
\begin{cases}
\dim P^-
&\quad \text{ if }
n\equiv 0\bmod 4,\\
\dim P^+
&\quad \text{ if }
n\equiv 2\bmod 4
\end{cases}
\label{eqcHdg}
\end{equation}
modulo $4$.

Let $v^{q-}$ be
the multiplicity
of the eigenvalue
$-1$ of the action
of the complex conjugate
on $P^q$
and let $(d^{q+},d^{q-})$
be the signature
of the quadratic vector
space $D^q$.
Then, 
we have
$e_q
-e_{q-2}
=v_q^-$
for $2\le q\le n$ even
and 
$e_0=v_0^-$.
Since $P^q$ is
$(-1)^{q/2}$-definite,
we have
$$v^{q-}
=
\begin{cases}
d^{q-}
&\quad \text{ if }
q\equiv 0\bmod 4,\\
d^{q+}
&\quad \text{ if }
q\equiv 2\bmod 4.
\end{cases}$$
If $n\equiv 0\bmod 4$,
we have
\begin{align*}
e_n^-
-
2e^-
\equiv
\sum_{q\le n,
\text{even}}
v_q^-
-2\sum_{q\le n,
q\equiv 2\ (4)}v_q^-
=&\
\sum_{q\le n,
q\equiv 0\ (4)}v_q^-
-\sum_{q\le n,
q\equiv 2\ (4)}v_q^-\\
=&\
\sum_{q\le n,
q\equiv 0\ (4)}d_q^-
-\sum_{q\le n,
q\equiv 2\ (4)}d_q^+
=
d^{\prime -}
-\dim P^-.
\end{align*}
If $n\equiv 2\bmod 4$,
we have
\begin{align*}
e_n^-
-
2e^-
\equiv
\sum_{q\le n,
\text{even}}
v_q^-
-2\sum_{q\le n,
q\equiv 0\ (4)}v_q^-
=&\
\sum_{q\le n,
q\equiv 2\ (4)}v_q^-
-\sum_{q\le n,
q\equiv 0\ (4)}v_q^-\\
=&\
\sum_{q\le n,
q\equiv 2\ (4)}d_q^+
-\sum_{q\le n,
q\equiv 0\ (4)}d_q^-
=
d^{\prime -}
-\dim P^+.
\end{align*}
Thus (\ref{eqcHdg})
is proved.

We have
$sw(H^n(X))=
(1+\{-1\})^{e_n^-}$
and
$hw(H^n_{dR}(X)')=
(1+\{-1\})^{d^{\prime -}}$.
Hence by (\ref{eqcHdg}),
we obtain
$$sw(H^n(X))
\cdot
(1+e^-\{-1,-1\})
=
hw(H^n_{dR}(X)')
\cdot
\begin{cases}
(1+\{-1\})^{-\dim P^-}
&\quad \text{ if }
n\equiv 0\bmod 4,\\
(1+\{-1\})^{-\dim P^+}
&\quad \text{ if }
n\equiv 2\bmod 4.
\end{cases}$$
Further, by
(\ref{eqLefr}),
we have
$$sw(H^n(X))
\cdot
(1+(e^-+\beta)\{-1,-1\})
=
hw(H^n_{dR}(X)')
\cdot
(1+\{-1\})^r.$$
Since
$e=e^-\{-1\}$,
the assertion is proved.
\end{proof}

\section{Stiefel-Whitney classes
of symmetric complexes}\label{sswcx}

In this section,
we define and
study the Stiefel-Whitney classes
of symmetric complexes.
A similar framework is studied
in \cite{Ba}. After recalling
some elementary constructions
on complexes,
we prepare
basic properties
on symmetric 
strict perfect complexes
in (\ref{dfscx})--(\ref{lmLag})
and we recall fundamental
properties of
the Stiefel-Whitney classes
of symmetric bundles
in (\ref{prEKV})--(\ref{corEKV}).
We define 
the class for
symmetric 
strict perfect complexes
in Definition \ref{dfw} and 
establish
fundamental properties in
Proposition \ref{prw}.
Finally,
we generalize the definition
to the derived category
$D_{\rm perf}(X)$
in Corollary \ref{corwD}.

\setcounter{thm}1

We will follow the sign
convention on complexes
in \cite{BBM}.
Let ${\mathcal  C}$ be 
an abelian category.
The mapping cone
${\rm Cone}(f)$
of a morphism 
$f\colon K\to L$
of complexes
of objects of ${\mathcal  C}$
is defined to be the simple
complex associated to
the double complex
$\cdots \to 0\to
K\to L\to 0\to \cdots$
where $L$ is put
on the first degree 0.
More concretely,
the $i$-th component is
$K^{i+1}\oplus L^i$
and the map is
given by the left multiplication
by the matrix
$\begin{pmatrix}
-d^{i+1}& 0 \\
f^{i+1}& d^i
\end{pmatrix}.$
Canonical maps
$L\to {\rm Cone}(f)
\to K[1]$
are defined 
as $L={\rm Cone}(0\to L)
\to {\rm Cone}(f)
\to {\rm Cone}(K\to 0)
=K[1]$.

Similarly,
the mapping fiber
${\rm Fib}(f)$
is defined to be the simple
complex associated to
the double complex
$\cdots \to 0\to
K\to L\to 0\to \cdots$
where $K$ is put
on the first degree 0.
The $i$-th component is
$K^i\oplus L^{i-1}$
and the map is
given by the left multiplication
by the matrix
$\begin{pmatrix}
d^i& 0 \\
f^i& -d^{i-1}
\end{pmatrix}.$
Canonical maps
$L[-1]\to {\rm Fib}(f)
\to K$
are defined 
as $L[-1]={\rm Fib}(0\to L)
\to {\rm Fib}(f)
\to {\rm Fib}(K\to 0)=K$.
We identify 
${\rm Fib}(f)[1]$
with 
${\rm Cone}(-f)$
by the identity.

Let $f\colon K\to L$
and $g\colon L\to M$
be morphisms
of complexes
such that
the composition
$g\circ f$
is homotope to $0$.
Then, a homotopy
$t$ connecting
$g\circ f$ to 0
defines a map $(g,t)\colon
{\rm Cone}(f)
\to M$.
Recall that a homotopy $t$
consists of 
maps
$t^i\colon K^i\to 
M^{i-1}$
satisfying
$g^i\circ f^i
=t^{i+1}\circ d^i
+d^{i-1}\circ t^i$.
The $i$-th component
of the map
$(g,t)$ is given by
$t^{i+1}\oplus g^i
\colon
K^{i+1}\oplus L^i\to M^i$.
The composition
$L\overset{\rm can}
\to {\rm Cone}(f)
\overset{(g,t)}
\to M$
is $g$.
A bijection
of the set of homotopies
connecting $g\circ f$ to
0 to the set of
maps of complexes
${\rm Cone}(f)
\to M$
such that the composition
with $L\to {\rm Cone}(f)$
is $g$
is defined 
by sending $t$ to $(g,t)$.
A homotopy $t$
also corresponds
to a map
$(f,t)\colon
K\to {\rm Fib}(g)$
such that the composition
with the canonical map
${\rm Fib}(g)\to L$
is equal to $f$.

For a homotopy
$t$ connecting
$g\circ f$ to 0,
we define a complex
\setcounter{equation}0
\begin{equation}
C=C(K\overset f
\to L\overset g
\to M)_t
\label{eqhtp}
\end{equation}
to be
${\rm Fib}((g,t)\colon
{\rm Cone}(f)\to M)=
{\rm Cone}((f,t)\colon
K\to {\rm Fib}(g))$.
For an integer $i$,
the $i$-th component
$C^i$ is given by
$C^i=K^{i+1}\oplus 
L^i\oplus M^{i-1}$
and the map
$d^i\colon C^i\to
C^{i+1}$
is given by the matrix
$\begin{pmatrix}
-d^{i+1}& 0 &0\\
f^{i+1}& d^i &0\\
t^{i+1}&g^i&-d^{i-1}
\end{pmatrix}.$
For a commutative diagram
$$\begin{CD}
K@>f>>L@>g>>M\\
@VaVV @VbVV @VcVV\\
K'@>{f'}>>L'@>{g'}>>M'
\end{CD}$$
of morphisms
of complexes
and homotopies $t$ and $t'$
connecting
$g\circ f$ and 
$g'\circ f'$ 
to 0 respectively
satisfying
$c\circ t=t'\circ a$,
the maps
$a,b$ and $c$ induce
a map
\begin{equation}
\begin{CD}
C=C(K\overset f
\to L\overset g
\to M)_t
@>>>
C'=C(K'\overset {f'}
\to L'\overset {g'}
\to M')_{t'}
\end{CD}
\label{eqCC'}
\end{equation}
of complexes.

Further we consider
a commutative diagram
$$\begin{CD}
K'@>{f'}>>L'@>{g'}>>M'\\
@V{a'}VV @V{b'}VV @V{c'}VV\\
K''@>{f''}>>L''@>{g''}>>M''
\end{CD}$$
of morphisms
of complexes
such that
$a'\circ a$,
$b'\circ b$ and
$c'\circ c$ are 0
and a homotopy $t''$
connecting
$g''\circ f''$ 
to 0 satisfying
$c'\circ t'=t''\circ a'$.
Then, we obtain 
maps
$C\to C'\to C''=
C(K''\overset{f''}\to
L''\overset{g''}\to
M'')_{t''}$
and the composition is 0.
The complex
$C^\circ=C(C\to C'\to C'')$
is the simple
complex associated
to the double complex
$\cdots \to0\to
C\to C'\to C''
\to 0\to\cdots$
where $C'$
is put on the first degree $0$.
We put
$K^\circ
=C(K\overset a\to 
K'\overset{a'}\to K''),
L^\circ
=C(L\overset b\to 
L'\overset{b'}\to L''),
M^\circ
=C(M\overset c\to 
M'\overset{c'}\to M'')$
and let
$f^\circ \colon
K^\circ
\to L^\circ$
and
$g^\circ \colon
L^\circ
\to M^\circ$
be the induced maps.
The homotopies
$t,t'$ and $t''$
induce a homotopy
$t^\circ$
connecting
$g^\circ
\circ f^\circ $ to $0$.
A canonical isomorphism
\begin{equation}
\begin{CD}
C^\circ=C(C\to C'\to C'')
@>>>
C(K^\circ
\overset{f^\circ}\to
L^\circ
\overset{g^\circ}\to
M^\circ)_{t^\circ}
\end{CD}
\label{eqCcirc}
\end{equation}
is defined by 
$1$ on
$L,K',L',M',L''$
and 
$-1$ on
$K,M,K'',M''$.

Let $D\colon 
{\mathcal  C}\to {\mathcal  C}$
be a contravariant
additive functor
and $c\colon {\rm id}
\to DD$
be a morphism of functors.
In practice,
the category
${\mathcal  C}$
will be that of
${\mathcal  O}_X$-modules
on a ringed space
$(X,{\mathcal  O}_X)$,
the functor $D$
will be defined by
$D{\mathcal  F}=
{\mathcal  H}om({\mathcal  F},{\mathcal  O}_X)$
and $c$ will be
defined by
the canonical morphism
${\mathcal  F}\to
DD{\mathcal  F}$.
For a complex $K$
of objects of ${\mathcal  C}$,
let $DK$
denote the dual complex.
The $i$-th component
$(DK)^i$ is equal to
$D(K^{-i})$
and $d^i\colon 
D(K^{-i})\to D(K^{-i-1})$
is given by
$(-1)^{i+1}D(d^{-i-1})$.
The $i$-th component
of the bidual complex $DDK$
is equal to $DD(K^i)$
and $d^i$
is given by $-DD(d^i)$.
The canonical map
$c_K\colon K\to DDK$ is defined
by $(-1)^i$-times the canonical map
$K^i\to DD(K^i)$.
The composition
$Dc_K\circ c_{DK}
\colon
DK\to DDDK
\to DK$
is the identity.

We define a canonical
isomorphism
$D{\rm Cone}(f)
\to 
{\rm Fib}(Df)$
by $(1,(-1)^{-i})$
on the $i$-th component
$D(L^{-i})\oplus
D(K^{-i+1})$.
In the case where
$L=0$,
this gives a canonical
isomorphism
$D(K[1])\to (DK)[-1]$.
Consequently,
for an even integer $n$
the canonical
isomorphism
$D(K[n])\to (DK)[-n]$
is defined
by the multiplication
by $(-1)^{n/2}$.
Similarly,
we define a canonical
isomorphism
$D{\rm Fib}(f)
\to 
{\rm Cone}(Df)$
by $(1,(-1)^{-i+1})$
on the $i$-th component
$D(K^{-i})\oplus
D(L^{-i-1})$.
The diagram
$$\begin{CD}
{\rm Cone}(f)
@>{{\rm Cone}(c_f)}>>
{\rm Cone}(DDf)
\\
@V{c_{{\rm Cone}(f)}}VV
@AAA\\
DD{\rm Cone}(f)
@<<<
D{\rm Fib}(Df)
\end{CD}$$
is commutative.
A similar diagram
obtained by
switching 
${\rm Cone}$
and
${\rm Fib}$
is also commutative.

The dual homotopy
$Dt$ connecting
$Df\circ Dg$ to 0
consists of
$(Dt)^i\colon D(M^{-i})
\to D(K^{1-i})$
defined by
$(Dt)^i=(-1)^iD(t^{1-i})$.
We define a canonical 
isomorphism
\begin{equation}
\begin{CD}
D(C(K\overset f
\to L\overset g
\to M)_t)
@>>>
C(DM\overset {Dg}
\to DL\overset {Df}
\to DK)_{Dt}
\end{CD}
\label{eqDC}
\end{equation}
to be the composition of
$$\begin{CD}
D({\rm Fib}((g,t)\colon
{\rm Cone}(f)\to M)
@>>>
{\rm Cone}(D(g,t)\colon
DM\to D{\rm Cone}(f))\\
@>>>
{\rm Cone}((Dg,Dt)\colon
DM\to {\rm Fib}(Df)).
\end{CD}
$$
It is defined by
$(-1)^i\oplus 1 \oplus
(-1)^{i+1}$
on the $i$-th component
$D(M^{-1-i})
\oplus D(L^{-i})
\oplus D(K^{1-i})$.
The diagram
\begin{equation}
\begin{CD}
C(K\overset f
\to L\overset g
\to M)_t
@>>>
C(DDK\overset {DDf}
\to DDL\overset {DDg}
\to DDM)_{DDt}
\\
@VVV @AAA\\
DDC(K\overset f
\to L\overset g
\to M)_t
@<<<
DC(DM\overset {Dg}
\to DL\overset {Df}
\to DK)_{Dt}
\end{CD}
\label{eqDDC}
\end{equation}
is commutative.

Let $(X,{\mathcal  O}_X)$
be a local ringed space
in the rest of this section.
Recall that a complex
$K=(K^i,d^i)$ of
${\mathcal  O}_X$-modules
is called a strict perfect complex
if $K^i$ is locally free of finite rank 
for each $i\in {\mathbb Z}$
and if $K^i$ are 0 except 
for finitely many $i\in {\mathbb Z}$.
A strictly perfect complex
is acyclic if and only
if it is locally homotope to $0$.
A map $f$
of strict perfect complexes
is a quasi-isomorphism
if and only if the mapping cone
${\rm Cone}(f)$
is acyclic.
If $K$ is
a strict perfect complex,
the canonical map
$c_K\colon K\to DDK$ is
an isomorphism.

\begin{df}\label{dfscx}
Let
$K$ be a strict perfect complex 
of ${\mathcal  O}_X$-modules.

{\rm 1.}
We say a morphism
$q:K\to DK$
of complexes
is symmetric 
if the composition
$Dq\circ c_K\colon
K\to DDK\to DK$ is 
equal to $q$.

{\rm 2.}
If a symmetric
morphism
$q:K\to DK$
is a
quasi-isomorphism,
we call the pair
$(K,q)$ 
a symmetric strict perfect complex.

{\rm 3.}
Let $f\colon
L\to DL$
be a symmetric morphism.
We say a homotopy
$t$ connecting
$f$ to $0$ is symmetric
if we have
$Dt\circ c_L=t$.
\end{df}

Let $(K,q)$ be
a symmetric strict
perfect complex.
Let $f\colon
L \to K$
be a morphism
of strict perfect
complexes
and let
$t$ be a symmetric homotopy
connecting 
the symmetric morphism
$Df\circ q\circ f\colon
L\to K\to DK\to DL$
to $0$.
We define a complex
\setcounter{equation}0
\begin{equation}
\begin{CD}
M=
M(L\overset f
\to K)_{q,t}
\end{CD}\label{eqM}
\end{equation}
to be $C(L\overset f
\to K
\overset{Df\circ q}
\to DL)_t$
{\rm (\ref{eqhtp})}
defined by the homotopy $t$.
We define a map
\begin{equation}
\begin{CD}
q_M\colon
M@>>> DM
\end{CD}
\label{eqqM}
\end{equation}
to be the composition
of the map 
$$\tilde q\colon
M=C(L\to K\to DL)_t
\to
C(DDL\to DK\to DL)_{Dt}
\quad
(\ref{eqCC'})$$
defined by the commutative
diagram
\begin{equation}
\begin{CD}
L@>f>> K@>{Df\circ q}>> DL\\
@V{c_L}VV @VqVV @|\\
DDL@>{Dq\circ DDf}>>
DK@>{Df}>>DL.
\end{CD}
\label{eqdqM}
\end{equation}
with the inverse of
the isomorphism
$$DM
=
DC(L\to K\to DL)_t
\to
C(DDL\to DK\to DL)_{Dt}
\quad
{\rm  (\ref{eqDC})}.$$

\begin{lm}\label{lmM}
Let $(K,q)$ be
a symmetric strict
perfect complex.
Let $f\colon
L \to K$
be a morphism
of strict perfect
complexes
and let
$t$ be a symmetric homotopy
connecting 
the symmetric morphism
$Df\circ q\circ f\colon
L\to K\to DK\to DL$
to $0$.

Then, the pair
$(M,q_M)$
is a symmetric 
strict perfect complex.
\end{lm}

\begin{proof}
The morphism
$q_M$ is a quasi-isomorphism
since the vertical arrows
in (\ref{eqdqM}) are
quasi-isomorphisms.
We show $q_M=Dq_M\circ c_M$.
By (\ref{eqDDC}),
the upper square of
the diagram
\setcounter{equation}0
\begin{equation}
\begin{CD}
M=
C(L\to K\to DL)_t
@>{\tilde c}>>
C(DDL\to DDK\to DDDL)_{DDt}
\\
@V{c_M}VV
@AAA\\
DDM=DDC(L\to K\to DL)_t
@<<<
DC(DDL\to DK\to DL)_{Dt}
\\
@.@V{D\tilde q}VV\\
@.
DM=
DC(L\to K\to DL)_t
\end{CD}
\label{eqCc}
\end{equation}
is commutative.
The composition
$DDM\to DM$
is $Dq_M$.
Let 
$$\widetilde{Dq}
\colon
C(DDL\to DDK\to DDDL)_{DDt}
\to
C(DDL\to DK\to DL)_{Dt}$$
be the map defined by
the dual diagram
\begin{equation}
\begin{CD}
DDL@>{DDf}>> DDK
@>{DDDf\circ DDq}>> DDDL\\
@| @V{Dq}VV @VV{Dc_L}V\\
DDL@>{Dq\circ DDf}>>
DK@>{Df}>>DL
\end{CD}
\label{eqDqM}
\end{equation}
of (\ref{eqqM}).
Then, the diagram
$$\begin{CD}
DC(DDL\to DK\to DL)_{Dt}
@>>>
C(DDL\to DDK\to DDDL)_{DDt}
\\
@V{D\tilde q}VV
@VV{\widetilde{Dq}}V\\
DM=
DC(L\to K\to DL)_t
@>>>
C(DDL\to DK\to DL)_{Dt}
\end{CD}$$
is commutative.
Since the composition
$DL\to DDDL\to DL$
is the identity,
we have $\tilde q
=\widetilde{Dq}\circ \tilde c$
where $\tilde c$
is the upper horizontal arrow
in (\ref{eqCc}).
Thus, we obtain
$q_M=Dq_M\circ c_M$.
\end{proof}

In the notation
$M(L\to K)_{q,t}$,
if the homotopy $t$ is 0,
we drop it from the notation
to make
$M(L\to K)_q$.

\begin{cor}\label{corK0}
Let $(K,q)$
be a symmetric 
strict perfect complex.
Let $K^{>0}\subset K$
denote the subcomplex
defined by
$(K^{>0})^i=K^i$ for $i>0$
and
$(K^{>0})^i=0$ for $i\le 0$
and define
a symmetric
strict perfect complex
$(K^\natural,q_{K^\natural})$
to be 
$M(K^{>0}\to K)_q$
defined in 
{\rm (\ref{eqM})}
and 
{\rm (\ref{eqqM})}.

Then, the complex
$K^\natural$ 
is acyclic except
at degree $0$
and the cohomology sheaf
${\mathcal  E}=
{\mathcal  H}^0(K^\natural)$
is locally free
of finite rank.
The map $q_{K^\natural}
\colon
K^\natural
\to DK^\natural$
induces
a non-degenerate
symmetric bilinear form
$q_{\mathcal  E}
\colon {\mathcal  E}\to
D{\mathcal  E}$.
\end{cor}

\begin{proof}
The complex 
$K^\natural$
is quasi-isomorphic
to the strict perfect complex
${\rm Cone}(K^{>0}\to (DK)^{\ge 0})$
supported on
degree $\ge 0$.
Since it is also
quasi-isomorphic
to the complex
${\rm Fib}(K^{\le 0}
\to D(K^{>0}))$,
it is acyclic
at degree $>0$.
Hence the assertion follows.
\end{proof}

By abuse of terminology,
we call a locally free
${\mathcal  O}_X$-module
of finite rank
a bundle on $X$.
For a locally free
${\mathcal  O}_X$-module
${\mathcal  E}$
of finite rank,
we call a locally free
${\mathcal  O}_X$-submodule
${\mathcal  F}$
a subbundle
if the quotient
${\mathcal  E}/{\mathcal  F}$
is locally free,
or equivalently,
if ${\mathcal  F}$
is locally a
direct summand of
${\mathcal  E}$.

\begin{df}
Let ${\mathcal  E}$
be a locally free
${\mathcal  O}_X$-module
of finite rank.

{\rm 1.}
We say an ${\mathcal  O}_X$-linear
map $q\colon
{\mathcal  E}\to
D{\mathcal  E}=
{\mathcal  H}om({\mathcal  E},{\mathcal  O}_X)$
is symmetric
if $q$ is equal to
the composition
$Dq\circ c_{\mathcal  E}
\colon 
{\mathcal  E}\to DD{\mathcal  E}
\to D{\mathcal  E}$.

{\rm 2.}
If a symmetric map $q\colon
{\mathcal  E}\to
D{\mathcal  E}$
is an isomorphism,
we say that $q$ is
non-degenerate and
we call the pair
$({\mathcal  E},q)$
a symmetric bundle.

{\rm 3.}
Let $({\mathcal  E},q)$
be a symmetric bundle.
We say a subbundle
${\mathcal  F}$
is isotropic
if the composition
${\mathcal  F}\to {\mathcal  E}
\to D{\mathcal  E}
\to D{\mathcal  F}$
is the $0$-map.
Further,
if the sequence
${\mathcal  F}\to {\mathcal  E}
\to D{\mathcal  F}$
is exact,
we say ${\mathcal  F}$
is Lagrangean.
\end{df}

If ${\mathcal  F}$
is an isotropic
subbundle
of a symmetric bundle
$({\mathcal  E},q)$,
the subquotient
\setcounter{equation}0
\begin{equation}
{\mathcal  E}'=
{\mathcal  H}({\mathcal  F}
\to {\mathcal  E}
\to D{\mathcal  F})=
{\rm Ker}({\mathcal  E}
\to D{\mathcal  F})/
{\mathcal  F}
\label{eqH}
\end{equation}
is locally free
and $q$
induces a
non-degenerate
symmetric bilinear form
$q'\colon
{\mathcal  E}'
\to D{\mathcal  E}'$.

\begin{pr}\label{prM}
Let $(K,q)$ be
a symmetric strict
perfect complex.
Let $f\colon
L \to K$
be a morphism
of strict perfect
complexes
and let
$t$ be a symmetric homotopy
connecting 
$Df\circ q\circ f$
to $0$.
Let $(M,q_M)$
be the symmetric
strict perfect complex
$M=
M(L\overset f
\to K)_{q,t}$
defined by
{\rm (\ref{eqM})} and
{\rm (\ref{eqqM})}.
Further we define $M^\natural=
M(M^{>0}\to M)_{q_M}$
as in Corollary 
{\rm \ref{corK0}}.

Assume that
$K$ is acyclic except
at degree $0$ and
that
the cohomology sheaf
${\mathcal  H}^0(K)={\mathcal  E}$
is locally free.
Then there exists
an isotropic subbundle
${\mathcal  F}$ of 
${\mathcal  H}^0(M^\natural)={\mathcal  E}'$
satisfying
the following properties:
\begin{itemize}
\item[{\rm (1)}]
There exists
an isomorphism
${\mathcal  H}({\mathcal  F}
\to {\mathcal  E}'
\to D{\mathcal  F})
\to {\mathcal  E}$
of symmetric bundles.
\item[{\rm (2)}]
There exists
an exact sequence
$0\to DL^1
\to {\mathcal  F}
\to D{\rm Ker}(K^1\to K^2)
\to 0$
of locally free
${\mathcal  O}_X$-modules.
\end{itemize}
\end{pr}

\begin{proof}
First, we prove
the case where
$L=0$.
In this case,
we have $M=K$
and $M^\natural
=C(K^{>0}\to K\to DK^{>0})$.
By the assumption
that
$K$ is acyclic except
at degree $0$,
the subcomplex
$K^{>0}$
is also acyclic
except at degree 1
and ${\mathcal  H}^1(
K^{>0})
={\rm Ker}(K^1\to K^2)$
is locally free of
finite rank.
Hence
$M^\natural$
is also acyclic
except at degree 0.
Further, the
cohomology sheaf
${\mathcal  E}'=
{\mathcal  H}^0(
M^\natural)$
is locally free
and 
has an increasing
3 step 
filtration $F^\bullet$
by subbundle
such that
${\rm Gr}_F^1
{\mathcal  E}'=
{\mathcal  H}^1(
K^{>0})$,
${\rm Gr}_F^0
{\mathcal  E}'
=
{\mathcal  E}=
{\mathcal  H}^0(K)$
and
${\rm Gr}_F^{-1}
{\mathcal  E}'
={\mathcal  H}^{-1}(
DK^{>0})=
D{\mathcal  H}^1(
K^{>0})$.
Since the composition
$(DK^{>0})[1]\to
M^\natural
=C(K^{>0}\to K\to DK^{>0})
\to K^{>0}[-1]$
is the zero map,
the subbundle
${\mathcal F}=
{\rm Gr}_F^{-1}
{\mathcal  E}'
={\mathcal  H}^{-1}(
DK^{>0})$
is isotropic
and the assertion (1) follows.
The assertion (2) follows from
${\mathcal  H}^1(
K^{>0})
={\rm Ker}(K^1\to K^2)$.

We prove the general case.
The complex
$L^\circ
=C(L^{>1}\to L\to D(L^{>1}))$
is acyclic except
at degree 1
and we have
${\mathcal  H}^1(
L^\circ)=L^1$.
By (\ref{eqCcirc}),
we have a
canonical isomorphism
$M^\natural
\to
C(L^\circ \to
K^\natural
\to DL^\circ)$.
Hence
${\mathcal  E}'=
{\mathcal  H}^0(
M^\natural)$
is locally free
and 
has an increasing
3 step 
filtration $F^\bullet$
by subbundle
such that
${\rm Gr}_F^1
{\mathcal  E}'=
{\mathcal  H}^1(
L^\circ)=L^1$,
${\rm Gr}_F^0
{\mathcal  E}'
=
{\mathcal  H}^0(K^\natural)$
and
${\rm Gr}_F^{-1}
{\mathcal  E}'
={\mathcal  H}^{-1}(
D(L^\circ))
=DL^1$.
Similarly as in the special case proved above,
the subbundle
${\rm Gr}_F^{-1}
{\mathcal  E}'
=DL^1$
is isotropic.
Further since it has been already proved
for $0\to K^\natural$,
there exist an isotropic subbundle
${\mathcal F}'=
{\mathcal  H}^1(
K^{>0})
={\rm Ker}(K^1\to K^2)$
of 
${\rm Gr}_F^0
{\mathcal  E}'
=
{\mathcal  H}^0(K^\natural)$
an isomorphism
${\mathcal H}
({\mathcal F}'
\to 
{\rm Gr}_F^0
{\mathcal  E}'
\to 
D{\mathcal F}')
\to {\mathcal E}$
of symmetric bundles.
Thus, the inverse image ${\mathcal F}
\subset
{\mathcal E}'$
of ${\mathcal F}'$
is isotropic and is an extension of
${\mathcal F}'={\rm Ker}(K^1\to K^2)$
by
${\rm Gr}_F^{-1}
{\mathcal  E}'
=DL^1$.
Thus the assertion follows.
\end{proof}

\begin{df}\label{dfsLag}
Let $(K,q)$ be
a symmetric strict
perfect complex.
Let $f\colon
L \to K$
be a morphism
of strict perfect
complexes.

If there exists
a symmetric homotopy $t$
connecting
$Df\circ q\circ f$
to $0$ such that
the complex
$M=M(L\overset f
\to K)_{q,t}$
{\rm (\ref{eqM})}
is acyclic,
we say $f\colon 
L\to K$
is Lagrangean.
\end{df}

The condition
$M=M(L\overset f
\to K)_{q,t}$
is acyclic
is equivalent
to that
the map
$(Df\circ q,t)\colon
{\rm Cone}(f)
\to DL$
is a quasi-isomorphism.

\begin{lm}\label{lmLag}
Let $(K,q)$ be
a symmetric strict
perfect complex.
Let $f\colon
L \to K$
be a morphism
of strict perfect
complexes
and let
$t$ be a symmetric homotopy
connecting 
$Df\circ q\circ f$
to $0$.
We put
$M=M(L\to K)_{q,t}$
and
$N={\rm Fib}(K\to DL)$.
Then the direct
sum $N\to K\oplus M$
of the maps defined by
$$
\begin{CD}
0@>>> K@>>> DL\\
@VVV @|@VVV\\ 
0@>>> K@>>> 0,
\end{CD}
\qquad
\begin{CD}
0@>>> K@>>> DL\\
@VVV @|@|\\ 
L@>>> K@>>> DL
\end{CD}$$
is Lagrangean
with respect to
the symmetric bilinear form
$q\oplus -q_M$
on $K\oplus M$.
\end{lm}

\begin{proof}
It suffices
to show
that the complex
$C(N\to K\oplus M\to DN)_0$
is acyclic.
Let $K^\circ$
denote the acyclic
complex
$C(K\to K\oplus K\to DK)_0$
where 
$K\to K\oplus K$
is the diagonal map
and 
$K\oplus K\to DK$
is the difference
$(Dq,-Dq)$.
By (\ref{eqCcirc}),
we have an isomorphism
$C(N\to K\oplus M\to DN)_0
\to 
C({\rm Fib}({\rm id}_L)
\to K^\circ
\to {\rm Cone}({\rm id}_{DL}))_0$
and
the assertion follows.
\end{proof}

In the following,
we assume that $X$
is a scheme over ${\mathbb Z}[\frac12]$.
For a symmetric bundle
$({\mathcal  E},q)$
on $X$,
the Stiefel-Whitney class
$$w({\mathcal  E})
=1+w_1({\mathcal  E})
+w_2({\mathcal  E})+\cdots$$
is defined in
$H^*(X,{\mathbb Z}/2{\mathbb Z})
=
\prod_{i=0}^{\infty}
H^*(X,{\mathbb Z}/2{\mathbb Z})$,
see \cite[$\S$1]{EKV}.
For a locally free
${\mathcal  O}_X$-module
${\mathcal  F}$
of rank $r$,
we put
\setcounter{equation}0
\begin{equation}
\bar c({\mathcal  F})
=
\sum_{i=0}^r
(1+\{-1\})^{r-i}
c_i({\mathcal  F})
\label{eqcF}
\end{equation}
in $H^*(X,{\mathbb Z}/
2{\mathbb Z})$.

The Stiefel-Whitney class is characterized by
the following properties:
\begin{itemize}
\item[{\rm (\ref{prEKV}.1)}]
For a morphism
$f\colon X\to Y$ of
schemes and 
a symmetric bundle
$({\mathcal  E},q)$
on $Y$,
we have
$f^*w({\mathcal  E})
=w(f^*{\mathcal  E})$.
\item[{\rm (\ref{prEKV}.2)}]
For the direct sum
$({\mathcal  E}\oplus
{\mathcal  E}',
q\oplus q')$
of symmetric bundles
$({\mathcal  E},q)$
and 
$({\mathcal  E}',q')$
on $X$,
we have
$w({\mathcal  E}\oplus
{\mathcal  E}')
=w({\mathcal  E})
\cdot w({\mathcal  E}')$.
\item[{\rm (\ref{prEKV}.3)}]
Assume ${\mathcal  E}$ is of rank $1$
and let
$w_1({\mathcal  E})$
be the class of ${\mathcal  E}$
as an element of
$H^1(X,{\mathbf O}(1))=
H^1(X,{\mathbb Z}/2{\mathbb Z}).$
Then,
we have
$w({\mathcal  E})=1+
w_1({\mathcal  E})$.
\end{itemize}

\begin{pr}\label{prEKV}
Let $({\mathcal  E},q)$ be a symmetric bundle.

{\rm 1.}
We have
\setcounter{equation}3
\begin{equation}
w({\mathcal  E},q)\cdot
w({\mathcal  E},-q)=
\bar c({\mathcal  E}).
\label{eqw-}
\end{equation}

{\rm 2.}
Let ${\mathcal  F}$
be an isotropic subbundle
of ${\mathcal  E}$
and regard
${\mathcal  E}'=
{\mathcal  H}(
{\mathcal  F}
\to
{\mathcal  E}
\to D{\mathcal  F})$
{\rm (\ref{eqH})}
as a symmetric bundle
by the symmetric bilinear
form $q'$ induced by $q$.
Then, we have
\begin{equation}
w({\mathcal  E})
=
w({\mathcal  E}')
\cdot
\bar c({\mathcal  F})
\end{equation}
\end{pr}

\begin{proof}
We define a map
${\mathcal  F}'=
{\rm Ker}({\mathcal  E}
\to D{\mathcal  F})
\to
{\mathcal  E}\oplus {\mathcal  E}'$
to be
the sum of the inclusion
and the projection.
Then, it is Lagrangean
with respect to
the symmetric bilinear form
$q\oplus-q'$.
Hence, by
(\ref{prEKV}.2)
and 
\cite[Proposition 5.5]{EKV},
we obtain
$w({\mathcal  E},q)\cdot
w({\mathcal  E}',-q')=
\bar c({\mathcal  F}')$.
By considering the
case where ${\mathcal  F}=0$,
we obtain the equality
(\ref{eqw-}).
Further, we obtain
$w({\mathcal  E},q)\cdot
w({\mathcal  E}',q')^{-1}=
\bar c({\mathcal  F}')\cdot
\bar c({\mathcal  E})^{-1}
=
\bar c({\mathcal  F})$.
\end{proof}

For a strict perfect
complex $K$,
we put
\begin{equation}
\bar c(K)
=
\prod_i
\bar c(K^i)^{(-1)^i}
\label{eqcbK}
\end{equation}
in $H^*(X,{\mathbb Z}/
2{\mathbb Z})$.
For a quasi-isomorphism 
$K_1\to K_2$,
we have
$\bar c(K_1)
=
\bar c(K_2)$.
For the dual complex $DK$,
we have
$\bar c(DK)
=
\bar c(K)$.

\begin{cor}\label{corEKV}
Let $(K,q)$ be
an acyclic
symmetric strict perfect
complex
and put ${\mathcal  E}=
{\mathcal  H}^0(K^\natural)$
for $K^\natural=
M(K^{>0}\to K)_q$.
Then, we have
\setcounter{equation}0
\begin{equation}
w({\mathcal  E})
\cdot \bar c(
K^{>0})=1.
\end{equation}
\end{cor}

\begin{proof}
We apply
Proposition \ref{prM}
to $L=0$.
Then ${\mathcal  F}=
D{\rm Ker}(K^1\to K^2)$
is a Lagrangean 
subbundle of ${\mathcal  E}$.
Hence, by Proposition
\ref{prEKV},
we obtain
$w({\mathcal  E})
=\bar c({\rm Ker}(
K^1\to
K^2))$.
Since $K$ is acyclic,
we have
$\bar c(K^{>0})=
\bar c({\rm Ker}(
K^1\to
K^2))^{-1}$.
\end{proof}

\begin{df}\label{dfw}
Let $(K,q)$
be a symmetric
strict perfect complex
and put ${\mathcal  E}=
{\mathcal  H}^0(K^\natural)$
for $K^\natural
=M(K^{>0}\to K)_q$.
Then, we define
the total Stiefel-Whitney class
$w(K)\in H^*(X,{\mathbb Z}/
2{\mathbb Z})$
by 
\setcounter{equation}0
\begin{equation}
w(K)
=w({\mathcal  E})\cdot
\bar c(K^{>0}).
\label{eqdfw}
\end{equation}
\end{df}

\begin{pr}\label{prw}
The Stiefel-Whitney classes
satisfy
the following properties:
\begin{itemize}
\item[{\rm (\ref{prw}.1)}]
For a morphism
$f\colon X\to Y$ of
schemes and 
a symmetric strict
perfect complex
$(K,q)$
on $Y$,
we have
$f^*w(K)
=w(f^*K)$.
\item[{\rm (\ref{prw}.2)}]
For the direct sum of
symmetric strict 
perfect complexes
$(K_1,q_1)$ and
$(K_2,q_2)$,
we have $w(K_1\oplus K_2)=
w(K_1)\cdot w(K_2)$.
\item[{\rm (\ref{prw}.3)}]
For a symmetric bundle
$({\mathcal  E},q)$ on $X$ and $K={\mathcal  E}[0]$,
we have $w(K)=w({\mathcal  E})$.
\item[{\rm (\ref{prw}.4)}]
For symmetric strict perfect complexes
$(K_1,q_1)$ and
$(K_2,q_2)$ of 
${\mathcal  O}_X$-modules
and 
a quasi-isomorphism
$f\colon K_1\to K_2$ 
such that the diagram
\setcounter{equation}6
\begin{equation}
\begin{CD}
K_1@>{f}>> K_2\\
@V{q_1}VV @VV{q_2}V\\
DK_1@<{Df}<< DK_2
\end{CD}
\label{eqsw1}
\end{equation}
is commutative up to homotopy,
we have $w(K_1)=w(K_2)$.
\item[{\rm (\ref{prw}.5)}]
For a symmetric strict perfect complex
$(K,q)$ and
a Lagrangean morphism 
$L\to K$ 
of strict perfect complexes,
we have $w(K)=\bar c(L)$.
\item[{\rm (\ref{prw}.6)}]
For a symmetric strict perfect complex
$(K,q)$,
we have $w(K,q)\cdot
w(K,-q)=\bar c(K)$.
\end{itemize}
\end{pr}
\begin{proof}
The properties
(\ref{prw}.1--3)
are clear from the definition. 
We show (\ref{prw}.6).
We put ${\mathcal  E}=H^0(K^\natural)$.
Then, by Proposition \ref{prEKV}.1,
we obtain
$$w(K,q)\cdot
w(K,-q)=
w({\mathcal  E},q_{\mathcal  E})\cdot
w({\mathcal  E},-q_{\mathcal  E})
\cdot
\bar c(K^{>0})^2
=
\bar c({\mathcal  E})
\cdot
\bar c(K^{>0})^2.$$
Since
$\bar c({\mathcal  E})
=
\bar c(K^\natural)
=
\bar c(K)
\cdot
\bar c(K^{>0})^{-2}$,
we obtain (\ref{prw}.6).
To show (\ref{prw}.5),
we first prove the following.

\begin{lm}\label{lmLagE}
Assume that $K$ is acyclic
except at degree $0$
and that
${\mathcal  E}={\mathcal  H}^0(K)$
is locally free of
finite rank. 
Let $f\colon
L\to K$
be a Lagrangean
morphism of strict
perfect complexes.
Then, we have
$w({\mathcal  E})=\bar c(L)$.
\end{lm}

\begin{proof}
Let $t$ be
a symmetric homotopy
connecting
$Df\circ q\circ f$
to 0.
We apply
Proposition \ref{prM}
to $L\to K$.
Then, by Proposition 
\ref{prEKV}.2,
we have
$w({\mathcal  E}')
=w({\mathcal  E})
\cdot
\bar c(L^1)
\cdot
\bar c(
{\rm Ker}(K^1\to K^2))$.
By the assumption
that $L$ is Lagrangean,
the complex $M=M(L\to K)_t$
is acyclic.
Hence, by
Corollary \ref{corEKV},
we have
$w({\mathcal  E}')
\cdot
\bar c(M^{>0})
=1$.
Since $K$ is acyclic except
at degree $0$,
we have
$\bar c(K^{>0})=
\bar c({\rm Ker}(K^1\to K^2))^{-1}$.
Thus, we obtain
$\bar c(K^{>0})
=w({\mathcal  E})
\cdot
\bar c(L^1)
\cdot
\bar c(M^{>0})$.
Since
$\bar c(M^{>0})=
\bar c(K^{>0})
\cdot
\bar c(L^{>1})^{-1}
\cdot
\bar c(D(L^{<0}))^{-1}$
and
$\bar c(L)=
\bar c(L^1)^{-1}
\cdot
\bar c(L^{>1})
\cdot
\bar c(D(L^{<0}))$,
the assertion follows.
\end{proof}

We prove 
(\ref{prw}.5).
Let $K^\natural$
denote the strict perfect
complex
$M(K^{>0}\to K)_q
=C(K^{>0}\to K
\to DK^{>0})_q.$
Then the map ${\mathcal  L'}
=C(0\to L
\to DK^{>0})
\to K^\natural$
induced by
$f\colon L\to K$ 
is Lagrangean.
By applying
Lemma \ref{lmLagE}
to $L'\to
K^\natural$,
we obtain
$w({\mathcal  E})
=\bar c(L')
=\bar c(L)
\cdot
\bar c(K^{>0})^{-1}$.
Thus,
we obtain
$w(K)=
w({\mathcal  E})
\cdot
\bar c(K^{>0})
=\bar c(L)$.

Finally, we deduce 
(\ref{prw}.4)
from (\ref{prw}.2),
(\ref{prw}.5) and
(\ref{prw}.6).
We consider
the direct sum
$(K_1\oplus K_2,q_1\oplus -q_2)$
as a symmetric strict perfect complex.
Let $t$ be a homotopy
connecting
$Df\circ q_2\circ f$ to $q_1$.
By replacing 
$t$ by the homotopy
$(t+Dt\circ c_{K_2})/2$,
we may assume $t$
is symmetric.
Then, the map
$({\rm id}_{K_1},f)\colon 
K_1\to K_1\oplus K_2$
is Lagrangean
with respect to
$-q_1\oplus q_2$.
Hence, by (\ref{prw}.2) and 
(\ref{prw}.5),
we obtain
$w(K_1,-q_1)\cdot
w(K_2,q_2)=
\bar c(K_2)$.
Thus,
we obtain
$w(K_1,q_1)=
w(K_2,q_2)$ by
(\ref{prw}.6).
\end{proof}

The properties
(\ref{prw}.2--5)
characterize the
Stiefel-Whitney classes.
We deduce
the uniqueness from
Lemma \ref{lmLag}
applied to
$L=K^{>0}$.
Since the symmetric strict
perfect complex
$K^\natural
=M(K^{>0}\to K)_q$
is quasi-isomorphic
to ${\mathcal  E}={\mathcal  H}^0(K^\natural)$, 
we obtain
$$w(K,q)\cdot w({\mathcal  E},-q_{\mathcal  E})
=w(K,q)\cdot w(K^\natural,
-q_{K^\natural})
=w(K\oplus K^\natural,
q\oplus -q_{K^\natural})
=\bar c(K^{>0}).$$

\begin{cor}\label{corw}
Let $(K,q)$ be
a symmetric strict
perfect complex.

{\rm 1.}
Let $f\colon
L \to K$
be a morphism
of strict perfect
complexes
and let
$t$ be a symmetric homotopy
connecting 
$Df\circ q\circ f$
to $0$.
We put
$M=M(L\to K)_{q,t}$.
Then, we have
$$w(K)=w(M)\cdot \bar c(L).$$

{\rm 2.}
Assume that
${\mathcal  H}^i(K)$
is locally free of
finite rank for 
every $i\in {\mathbb Z}$
and regard
${\mathcal  H}^0(K)$
as a symmetric bundle
by the symmetric form
induced by $q$.
Then, we have
$$w(K)=w({\mathcal  H}^0(K))\cdot 
\prod_{i<0}
\bar c({\mathcal  H}^i(K))
^{(-1)^i}.$$
\end{cor}

\begin{proof}
1.
We put
$N={\rm Fib}(K\to DL)$
as in Lemma \ref{lmLag}.
Then, by (\ref{prw}.2),
(\ref{prw}.5),
and by Lemma \ref{lmLag},
we have
$w(K,q)
\cdot w(M,-q_M)
=\bar c(N)$.
Further by (\ref{prw}.6),
we obtain
$w(K,q)
\cdot
\bar c(M)
=w(M,q_M)
\cdot
\bar c(N)$.
By 
$\bar c(M)
=
\bar c(L)^{-1}
\cdot
\bar c(N)$,
the assertion follows.

2.
We define a subcomplex
$L\subset K$
by
$L^i=K^i$ if $i<0$,
$L^0={\rm Im}
(K^{-1}\to K^0)$
and
$L^i=0$ if $i>0$.
Then, the inclusion
$i\colon L\to K$
satisfies
$Di\circ q\circ i=0$.
The complex $M=
M(L\to K)_{q,0}$
is acyclic except at
degree 0
and we have
${\mathcal  H}^0(M)=
{\mathcal  H}^0(K)$.
Thus,
the assertion follows
from 1.
\end{proof}

In the rest of this section,
we assume that
the scheme $X$ is 
divisorial
\cite[Definition 2.2.5]{SGA6}
and
either separated
or noetherian.
Recall from
\cite[Proposition 2.2.9 b)]{SGA6}
that the natural functor
$K(X)\to
D_{\rm perf}(X)$
from the homotopy
category
of strict perfect complexes
of ${\mathcal  O}_X$-modules
to the derived category
of perfect complexes
of ${\mathcal  O}_X$-modules
induces an equivalence of
categories from the
quotient category
divided by quasi-isomorphisms.
For an object
$K$ of 
$D_{\rm perf}(X)$
and an isomorphism
$q\colon K\to DK$
satisfying
$q=Dq\circ c_K$,
we call the pair
$(K,q)$
a symmetric perfect
complex on $X$.

\begin{lm}\label{lmDK}
Let $X$ be a divisorial scheme
over ${\mathbb Z}[\frac12]$
either separated
or noetherian.
Let $K$ be an object of
the derived category
$D_{\rm perf}(X)$
of perfect complexes
of ${\mathcal  O}_X$-modules
and $q\colon K\to DK$
be an isomorphism of 
$D_{\rm perf}(X)$
satisfying $q=Dq\circ c_K$.

{\rm 1.}
There exist a 
symmetric strict perfect complex
$(K_0,q_0)$ and 
a quasi-isomorphism
$f_0\colon K_0\to K$ of 
${\mathcal  O}_X$-modules
such that the diagram
\setcounter{equation}0
\begin{equation}
\begin{CD}
K_0@>{f_0}>> K\\
@V{q_0}VV @VVqV\\
DK_0@<{Df_0}<< DK
\end{CD}
\label{eqKK0}
\end{equation}
in $D_{\rm perf}(X)$
is commutative.

{\rm 2.}
Let $(K_1,q_1)$ be another
symmetric strict perfect complex
and $f_1\colon K_1\to K$
be a quasi-isomorphism
such that the diagram
{\rm (\ref{eqKK0})}
with suffix $0$ replaced by $1$
is commutative.
Then, there exist
a strict perfect complex $K_2$
and quasi-isomorphisms
$g_0\colon K_2\to K_0$
and 
$g_1\colon K_2\to K_1$
such that the diagram
\begin{equation}
\begin{CD}
K_2@>{g_0}>> K_0@>{q_0}>>DK_0\\
@V{g_1}VV @. @VV{Dg_0}V\\
K_1@>{q_1}>>DK_1@>{Dg_1}>> DK_2
\end{CD}
\label{eqKK012}
\end{equation}
in commutative
up to homotopy.
\end{lm}

\begin{proof}
1.
There exist
strict perfect complexes
$K_1$ and $K_2$,
quasi-isomorphisms
$f_1\colon K_1\to K$ and
$f_2\colon K_2\to K$ 
and a morphism
$q_1\colon K_1\to DK_2$
such that the diagram
\begin{equation}
\begin{CD}
K_1@>{f_1}>> K\\
@V{q_1}VV @VVqV\\
DK_2@<{Df_2}<< DK
\end{CD}
\label{eqKK12}
\end{equation}
in $D_{\rm perf}(X)$
is commutative.
Further, there exist
a strict perfect complex
$K_0$ and
quasi-isomorphisms
$g_1\colon K_0\to K_1$ and
$g_2\colon K_0\to K_2$ 
such that the diagram
\begin{equation}
\begin{CD}
K_0@>{g_1}>> K_1\\
@V{g_2}VV @VV{f_1}V\\
K_2@>{f_2}>> K
\end{CD}
\label{eq012}
\end{equation}
in $D_{\rm perf}(X)$
is commutative.
The composition 
$q'_0=Dg_2\circ q_1\circ g_1
\colon K_0\to DK_0$
makes the diagram
(\ref{eqKK0})
commutative and hence is
a quasi-isomorphism.
Since it is homotope to
$Dq'_0\circ c_{K_0}$,
the map
$q_0=(q'_0
+Dq'_0\circ c_{K_0})/2
\colon K_0\to DK_0$
is symmetric and
is homotope to $q'_0$.

2.
We consider a 
commutative diagram (\ref{eq012})
with suffix 1, 2, 0
replaced by 0, 1, 2.
Then, we obtain a
diagram (\ref{eqKK012})
commutative up to homotopy.
\end{proof}

\begin{cor}\label{corwD}
Let $X$ and $(K,q)$ be
as in Lemma {\rm \ref{lmDK}}.
Then, for 
a symmetric strict perfect complex
$(K_1,q_1)$ and 
an isomorphism 
$f_1\colon K_1\to K$ 
as in Lemma {\rm \ref{lmDK}.1},
the Stiefel-Whitney class
$w(K_1)$
is independent of $(K_1,q_1)$.
\end{cor}

\begin{proof}
Let $(K_2,q_2)$ be
another symmetric strict perfect complex
and 
$f_2\colon K_2\to K$ 
be an isomorphism 
as in Lemma {\rm \ref{lmDK}.1}.
We consider
a strict perfect complex $K_0$,
quasi-isomorphisms
$g_1\colon K_0\to K_1$
and
$g_2\colon K_0\to K_2$
and a symmetric homotopy $t$
as in Lemma {\rm \ref{lmDK}.2}.
Then the composition
$q_0=Dg_1\circ q_1\circ g_1\colon
K_0\to DK_0$
defines a 
symmetric strict perfect complex.
By (\ref{prw}.1),
we obtain
$w(K_1)=w(K_0)=w(K_2)$.
\end{proof}

For a symmetric perfect
complex
$(K,q)$ on $X$,
the Stiefel-Whitney
class $w(K)$ is defined
as $w(K_0)$
by taking a quasi-isomorphism
$f_0\colon K_0\to K$
as in Lemma \ref{lmDK}.1.
It is well-defined
by Corollary \ref{corwD}.

We give a slight generalization.
Let $K$
be a perfect complex and
$n$ be an even integer.
We say that
an isomorphism
$q\colon K\to
DK[-2n]$ 
in the derived category
$D_{\rm perf}(X)$
is symmetric if
$q$ is equal to
the composition
\setcounter{equation}0
\begin{equation}
\begin{CD}
K@>{c_K}>> 
DDK
@>{\rm can}>>
D(DK[-2n])[-2n]
@>{Dq[-2n]}>>
DK[-2n].
\end{CD}
\label{eqKn}
\end{equation}
Let
$q\colon K\to
DK[-2n]$
be a symmetric isomorphism.
We put $K'=K[n]$
and define $q'\colon
K'\to DK'$
to be the composition
$$
\begin{CD}
K'=K[n]@>{q[n]}>> 
DK[-n]
@>{\rm can}>>
D(K[n])=DK'.
\end{CD}$$
Then, $(K',q')$
is a symmetric perfect complex
and the Stiefel-Whitney class
$w(K')$
is defined.

\begin{cor}\label{corwn}
Let $K$
be a perfect complex,
$n$ be an even integer
and $q\colon K\to
DK[-2n]$ be a
symmetric isomorphism
in the derived category
$D_{\rm perf}(X)$.
Let $(K[n],q')$ be
the symmetric
perfect complex
defined above.

Assume that
${\mathcal  H}^i(K)$
is locally free of
finite rank for 
every $i\in {\mathbb Z}$
and regard
${\mathcal  E}={\mathcal  H}^n(K)$
as a symmetric bundle
by the symmetric form
$q_{\mathcal  E}$ induced by $q$.
Then, we have
$$w(K)=w({\mathcal  E},
(-1)^{n/2}q_{\mathcal  E})\cdot 
\prod_{i<0}
\bar c({\mathcal  H}^i(K))
^{(-1)^i}.$$
\end{cor}

\begin{proof}
For $K'=K[n]$,
we have
${\mathcal  H}^0(K')
={\mathcal  H}^n(K)
={\mathcal  E}$.
By the sign convention
on the canonical isomorphism
$D(K[n])\to DK[-n]$,
the symmetric bilinear
form on ${\mathcal  E}$
induced by $q'$
is equal to
$(-1)^{n/2}q_{\mathcal  E}$.
Hence the assertion
follows from Corollary \ref{corw}.2.
\end{proof}

\section{Families}\label{sfam}

We generalize definitions
in Section \ref{scn}.
Let $S$ be a normal scheme
over ${\mathbb Z}[\frac 1{2\ell}]$
and $L$ be a finite extension of
${\mathbb Q}_\ell$.
We say that a smooth $L$-sheaf $V$ on $S$
is orthogonal if it is endowed with a
non-degenerate symmetric bilinear form $b:V\otimes V\to L$.
We define the Stiefel-Whitney class
$$sw_2(V)\in H^2(S,{\mathbb Z}/2{\mathbb Z})$$
in \'etale cohomology as follows.
It suffices to consider
the case where $S$ is connected.
In this case,
an orthogonal $L$-sheaf
$(V,b)$ is defined by an orthogonal
representation 
$\rho_{\bar x}:\pi_1(S,\bar x)\to O(V_{\bar x},b_{\bar x})$
for a base point $\bar x$.
Then we define 
the Stiefel-Whitney class
$sw_2(V,b)\in H^2(S,{\mathbb Z}/2{\mathbb Z})$
to be the image of
$sw_2(\rho_{\bar x})\in
H^2(\pi_1(S,{\bar x}),{\mathbb Z}/2{\mathbb Z})$
by the canonical map
$H^2(\pi_1(S,{\bar x}),{\mathbb Z}/2{\mathbb Z})
\to
H^2(S,{\mathbb Z}/2{\mathbb Z})$.
It is independent of the choice of base point $\bar x$.
When there is no fear of confusion,
we drop $b$ and write simply $sw_2(V)$.
Similarly as (\ref{eqdfswg}),
for a graded smooth $L$-sheaf 
$V^\bullet$
on $S$,
we define the Stiefel-Whitney
classes
$sw_1(V^\bullet)
\in H^1(S,
{\mathbb Z}/2{\mathbb Z})$
and
$sw_2(V^\bullet)
\in 
H^2(S,{\mathbb Z}/2{\mathbb Z})$.
 
The definition of the Stiefel-Whitney class
commutes with base change.
When $S$ is ${\rm Spec}\ K$ for a field $K$
and the smooth orthogonal $L$-sheaf $V$ on
$S$ is defined by an
orthogonal $L$-representation 
$\rho$ of $\text{Gal}(\bar K/K)$,
we have $sw_2(V)=sw_2(\rho)$
in $H^2(S,{\mathbb Z}/2{\mathbb Z})=H^2(\text{Gal}(\bar K/K),{\mathbb Z}/2{\mathbb Z})$.

The first Stiefel-Whitney class
$sw_1(V,b)\in H^1(S,{\mathbb Z}/2{\mathbb Z})$
of an orthogonal smooth $L$-sheaf $V$ on $S$
is similarly defined as follows.
We may assume $S$ is connected
and consider the corresponding representation
$\rho:\pi_1(S,\bar x)\to O(V)$.
Then, the class
$sw_1(V,b)$
is the homomorphism
$\det\rho:\pi_1(S,\bar x)\to \{\pm1\}\simeq
{\mathbb Z}/2{\mathbb Z}$
regarded as an element of
$H^1(S,{\mathbb Z}/2{\mathbb Z})
=
Hom(\pi_1(S,\bar x),\{\pm1\})$.

Let $S$ be a connected
normal scheme
over ${\mathbb Z}
[\frac 1{2\ell}]$.
Let $f:X\to S$ be a proper smooth morphism
of relative even dimension $n$.
Let $V^i$
denote the smoooth 
${\mathbb Q}_\ell$-sheaf
$R^{i+n}f_*
{\mathbb Q}_\ell(\frac n2)$
on $X$.
The cup-product 
defines a non-degenerate 
$(-1)^i$-symmetric bilinear form
$V^i\times V^{-i}
\to
R^{2n}f_*
{\mathbb Q}_\ell(n)
\to {\mathbb Q}_\ell$
by sending $(x,y)$ to
$\text{Tr}(x\cup y)$.
We define 
the second
Stiefel-Whitney classes
by
\setcounter{equation}0
\begin{align}
&sw_2(H^n_\ell(X/S))
=
sw_2(V^0),
\label{dfswXS}\\
&
sw_2(H^\bullet_\ell(X/S))
=
sw_2(V^\bullet)
=
sw_2(V^0)
+\sum_{q<0}
\bar c_1(V^q).
\nonumber
\end{align}
as elements
in $H^2(S,{\mathbb Z}/2{\mathbb Z})$,
similarly as in (\ref{eqdfswg}).

For an integer $q\ge 0$,
we consider the character
$e_q\colon
\pi_1(S)^{\rm ab}
\to
\{\pm1\}$
(\ref{eqeq})
as an element
of $H^1(S,{\mathbb Z}/2{\mathbb Z})$
and let $b_{{\acute et},q}$
be the rank of the
smooth ${\mathbb Q}_\ell$-sheaf
$R^qf_*
{\mathbb Q}_\ell$.
We put
$e=\sum_{q<n} e_q$
as in Conjecture \ref{cn}.
We also put
$\beta=\displaystyle{
\frac12\sum_{q<n}
(-1)^q(n-q)
b_{\text{\it \'et},q}}$.
By the definition
(\ref{eqdfswg})
and Lemma \ref{lmcchi}.1,
they are related
by the equality
\begin{equation}
sw_2(H^\bullet_\ell(X/S))
=
sw_2(H^n_\ell(X/S))
+\{e,-1\}
+\beta\cdot c_\ell.
\label{dfswXS2}
\end{equation}
The first Stiefel-Whitney class
$sw_1(H^n_\ell(X/S))
=sw_1(H^n_\bullet(X/S))
\in H^1(S,{\mathbb Z}/2{\mathbb Z})$
is defined
as
$\det R^nf_*
{\mathbb Q}_\ell(\frac n2)$.

We define the Hasse-Witt class
for the de Rham cohomology.
Let $S$ be a divisorial
scheme 
\cite[Definition 2.2.5]{SGA6}
over 
${\mathbb Z}[\frac1{2\ell}]$
either separated
or noetherian
and let $f:X\to S$ be a proper
smooth morphism
of relative even dimension $n$.
We put
${\mathcal  H}^\bullet_{dR}(X/S)
=Rf_*\Omega^\bullet_{X/S}$.
The product
$\Omega^\bullet_{X/S}
\otimes
\Omega^\bullet_{X/S}
\to
\Omega^\bullet_{X/S}$
and the trace map
$R^nf_*\Omega^n_{X/S}
\to {\mathcal  O}_S$
induce a symmetric isomorphism
${\mathcal  H}^\bullet_{dR}(X/S)
\to D{\mathcal  H}^\bullet_{dR}(X/S)
[-2n]$
by Poincar\'e duality.
Hence, the Hasse-Witt
classes
$hw_i({\mathcal  H}^\bullet_{dR}(X/S))
\in H^i(S,{\mathbb Z}/2
{\mathbb Z})$
are defined as
at the end of Section \ref{sswcx}.

\begin{lm}\label{lmD}
Assume that
${\mathcal H}^q_{dR}(X/S)
=R^qf_*\Omega^\bullet_{X/S}$ are locally
free for all $q
\in {\mathbb Z}$
and we put
$r=\sum_{q<n}(-1)^q 
{\rm rank} {\mathcal H}^q_{dR}(X/S)$.
Then ${\mathcal H}^n_{dR}(X/S)$
is a symmetric bundle
and
we have
\setcounter{equation}0
\begin{align}
&hw_2(H^\bullet_{dR}(X/S))\\
=&\
hw_2({\mathcal H}^n_{dR}(X/S))
\nonumber\\
&\qquad\qquad+
\begin{cases}
r\{d_X,-1\}+\binom r2\{-1,-1\}
& \text{if }n\equiv 0 \bmod 4,
\\
(r+b_{dR,n}-1)\{d_X,-1\}+
\binom {r+b_{dR,n}}2\{-1,-1\}
& \text{if }n\equiv 2 \bmod 4
\end{cases}
\nonumber\\
&\qquad\qquad
+
\sum_{0\le q<n}
c_1({\mathcal H}^q_{dR}(X/S))
\nonumber
\end{align}
in $H^2(S,{\mathbb Z}/2{\mathbb Z})$
and
$$hw_1(H^\bullet_{dR}(X/S))
=
hw_1({\mathcal H}^n_{dR}(X/S))
+
\frac n2\cdot 
b_{dR,n}\cdot \{-1\}
$$
in $H^1(S,{\mathbb Z}/2{\mathbb Z})$.
\end{lm}

The condition of Lemma
\ref{lmD}
is satisfied if
$S$ is the spectrum of
a field or
$S$ is a reduced scheme over ${\mathbb Q}$.

\begin{proof}
Similarly as Corollary \ref{cor'},
it follows from
Corollary \ref{corwn}
and the definition of $\bar c({\mathcal F})$
in (\ref{eqcF}).
\end{proof}

For the relation between
$sw_2(H^\bullet_\ell(X/S))$
and $hw_2(H^\bullet_{dR}(X/S))$,
we state the following conjecture.

\begin{cn}\label{cn2}
Let $S$ be a normal
divisorial scheme over 
${\mathbb Z}[\frac1{2\ell}]$
either separated
or noetherian
and let $f:X\to S$ be a proper
smooth morphism
of relative even dimension $n$.
We put
$$\eta=
\sum_{q<\frac n2}
(-1)^q\left(\frac n2-q\right)
{\rm rank}Rf_*\Omega^q_{X/S}$$
as in {\rm(\ref{eqr})}
and define $c_2,c_\ell
\in H^2(S,{\mathbb Z}/
2{\mathbb Z})$
as in {\rm(\ref{eqcl})}.

Then, we have an equality
\setcounter{equation}0
\begin{equation}
sw_2(H^\bullet_\ell(X/S))
= 
hw_2(H^\bullet_{dR}(X/S))
+\{2,hw_1(H^\bullet_{dR}(X/S))\}+
\eta\cdot (c_\ell-c_2).
\label{eqcncxS}
\end{equation}
in $H^2(S,{\mathbb Z}/2{\mathbb Z})$.
\end{cn}

If the condition
of Lemma \ref{lmD}
is satisfied,
the equality
(\ref{eqcncxS})
is equivalent to
\begin{align}
&\ sw_2(H^n_\ell(X/S))+ \{e,-1\}
+\beta\cdot c_\ell
\label{eqcnS}
\\
=&\ hw_2(H^n_{dR}(X/S))
\nonumber \\
&\qquad \qquad +
\begin{cases}\displaystyle{
r\{d_X,-1\}+\binom r2\{-1,-1\}}
& \text{if }n\equiv 0 \bmod 4,
\nonumber\\\displaystyle{
(r+b_{dR,n}-1)\{d_X,-1\}+
\binom {r+b_{dR,n}}2\{-1,-1\}}
& \text{if }n\equiv 2 \bmod 4
\end{cases}
\nonumber\\
&\qquad \qquad + \{2,d_X\}+
\eta\cdot(c_\ell-c_2)
+
\sum_{0\le q<n}
c_1({\mathcal H}^q_{dR}(X/S))
.\nonumber
\end{align}
Thus 
Conjecture \ref{cn2}
is a generalization
of Conjecture \ref{cn}.

The following weak evidence
on Conjecture \ref{cn2}
will be used in the
proof of the assertion 5
of Theorem \ref{thm1}
in the final section.

\begin{lm}\label{lmlS}
Let $S$ be a smooth scheme
over ${\mathbb Z}[\frac 12]$
and $f\colon X\to S$
be a proper smooth morphism
of even relative dimension $n$.
Let $\ell$ be a prime
number.
Assume that
the following condition is satisfied:
\begin{itemize}
\item[$(P)$]
For every finite
unramified extension
$K$ of ${\mathbb Q}_\ell$
and every
morphism
$f\colon
{\rm Spec}\ {\mathcal  O}_K\to
S$,
the pull-back
$X\times_SK$
by the restriction
$f|_K\colon
{\rm Spec}\ K\to S$
satisfies Conjecture {\rm \ref{cn}}.
\end{itemize}

Then,
the difference $\delta
\in H^2(S[\frac 1\ell],
{\mathbb Z}/2{\mathbb Z})$
of the both sides of
{\rm (\ref{eqcncxS})}
for $X[\frac 1\ell]\to S[\frac 1\ell]$
is in the image of
the restriction map
$H^2(S,
{\mathbb Z}/2{\mathbb Z})
\to
H^2(S[\frac 1\ell],
{\mathbb Z}/2{\mathbb Z})$.
\end{lm}

\begin{proof}
Since $S$ is assumed smooth
over ${\mathbb Z}[\frac 12]$,
we have an exact sequence
$$\begin{CD}
H^2(S,
{\mathbb Z}/2{\mathbb Z})
@>>>
H^2(S[\frac 1\ell],
{\mathbb Z}/2{\mathbb Z})
@>{\partial}>>
H^1(S_{{\mathbb F}_\ell},
{\mathbb Z}/2{\mathbb Z})
\end{CD}$$
by local acyclicity of
smooth morphism.
By a generalization of 
the Chebotarev density theorem 
\cite[Theorem 7]{ZL},
\cite[Theorem 9.11]{STai},
it suffices to show that
the image of $\delta$
by the composition
$$\begin{CD}
H^2(S[\frac 1\ell],
{\mathbb Z}/2{\mathbb Z})
@>{\partial}>>
H^1(S_{{\mathbb F}_\ell},
{\mathbb Z}/2{\mathbb Z})
\to
H^1(s,
{\mathbb Z}/2{\mathbb Z})
\end{CD}$$
is 0 for
every closed point $s$
of $S_{{\mathbb F}_\ell}$.

Let $K$
be an 
unramified extension
of ${\mathbb Q}_\ell$
with residue field
$\kappa(s)$.
Since $S$ is smooth,
the closed immersion
$s\to S$
is lifted to a morphism
$f\colon
{\rm Spec}\ {\mathcal  O}_K\to
S$.
Then, we have
a commutative diagram
$$\begin{CD}
H^2(S[\frac 1\ell],
{\mathbb Z}/2{\mathbb Z})
@>{\partial}>>
H^1(S_{{\mathbb F}_\ell},
{\mathbb Z}/2{\mathbb Z})
\\
@V{f|_K^*}VV @VVV
\\
H^2(K,
{\mathbb Z}/2{\mathbb Z})
@>{\partial}>>
H^1(s,
{\mathbb Z}/2{\mathbb Z}).
\end{CD}
$$
By the assumption $(P)$,
the pull-back
$X_K$ satisfies 
Conjecture \ref{cn}.
By the remark preceeding
Lemma \ref{lmlS},
it means
$f|_K^*(\delta)=0$
and the assertion follows.
\end{proof}

If $n=0$,
Conjecture \ref{cn2}
is nothing but
the following.

\begin{thm}[{{\rm \cite[Theorem 2.3]{EKV}}}]
For a finite \'etale morphism $f:X\to S$,
we have
$$sw_2(f_*{\mathbb Q})
=
hw_2(f_*{\mathcal  O}_X,{\rm Tr}_{X/S}(x^2))
+
\{2,d_X\}.$$
\end{thm}

\section{Transcendental argument}\label{strs}

In this section,
we prove the assertion
4 of Theorem \ref{thm1},
by a transcendental argument.
Since it will be reduced
to proving Conjecture \ref{cn2}
for schemes over ${\mathbb C}$,
we study them first.

We introduce some terminology.
Let $S$ be a connected
normal scheme 
of finite type over
${\mathbb C}$
and
let $S^{\rm an}$ denote the 
associated analytic space.
We say a local system $V$ of 
${\mathbb C}$-vector spaces 
on $S^{\rm an}$
is orthogonal if it is equipped with a
non-degenerate symmetric bilinear form 
$b:V\otimes V\to {\mathbb C}$.
An orthogonal local system
$V$ corresponds
to an orthogonal representation
$\rho_{\bar x}:\pi_1(S^{\rm an},\bar x)\to 
O(V_{\bar x},b_{\bar x})$
of the fundamental group
for a geometric point $\bar x$.
The Stiefel-Whitney class
$sw_2(V,b)\in H^2(S^{\rm an},{\mathbb Z}/2{\mathbb Z})$
is defined as
the image of
the Stiefel-Whitney class
$sw_2(\rho_{\bar x})\in 
H^2(\pi_1(S^{\rm an},\bar x),
{\mathbb Z}/2{\mathbb Z})$
by the canonical map
$H^2(\pi_1(S^{\rm an},\bar x),
{\mathbb Z}/2{\mathbb Z})
\to
H^2(S^{\rm an},{\mathbb Z}/2{\mathbb Z})$.
It is independent of 
the choice of base point $\bar x$.

A locally free ${\mathcal  O}_{S^{\rm an}}$-module $D$
is called quadratic if
it is equipped with a
non-degenerate symmetric bilinear form
$b:D\otimes D\to {\mathcal  O}_{S^{\rm an}}$.
For a quadratic locally free ${\mathcal  O}_{S^{\rm an}}$-module
$(D,b)$,
the Hasse-Witt class
$hw_2(D,b)\in H^2(S^{\rm an},{\mathbb Z}/2{\mathbb Z})$
is defined, see for example
\cite[\S1]{EKV}.
If $D$ is of rank $r$,
it is the image of
the class of
$D$ in $H^1(S^{\rm an},{\bf O}(r)^{\rm an})$
by the boundary map
$H^1(S^{\rm an},{\bf O}(r)^{\rm an})\to
H^2(S^{\rm an},{\mathbb Z}/2{\mathbb Z})$

\begin{lm}\label{lm1a}
Let $S$ be a normal scheme of finite type over ${\mathbb C}$.
Let $V$ be an orthogonal 
local system of ${\mathbb C}$-vector
spaces on
$S^{\rm an}$
and $D=V\otimes_{\mathbb C}
{\mathcal  O}_{S^{\rm an}}$ 
be the corresponding
quadratic locally free 
${\mathcal  O}_{S^{\rm an}}$-module.
Then we have
$sw_2(V)=hw_2(D)$ in 
$H^2(S^{\rm an},{\mathbb Z}/2 {\mathbb Z})$.
\end{lm}

\begin{proof}
We may assume $S$ is connected
and let $\pi_1(S^{\rm an},\bar b)$ be the
topological fundamental group
defined by a base point $\bar b$.
The isomorphism class
of the orthogonal local system $V$
is defined
as an element of
$H^1(S^{\rm an},O(r,{\mathbb C}))=
H^1(\pi_1(S^{\rm an},\bar b),O(r,{\mathbb C}))$.
By the definition of the Stiefel-Whitney class
$sw_2(V)$ and the functoriality
of the map from group cohomology
to singular cohomology,
it is equal to its image by the boundary map
$H^1(S^{\rm an},O(r,{\mathbb C})) \to H^2(S^{\rm an},{\mathbb Z}/2{\mathbb Z})$
defined by the central extension
$$\begin{CD}
1@>>>{\mathbb Z}/2{\mathbb Z}@>>> \tilde O(r,{\mathbb C})@>>> O(r,{\mathbb C})@>>> 1.
\end{CD}$$
Hence the assertion follows from 
the commutative diagram
$$\begin{CD}
1@>>>{\mathbb Z}/2{\mathbb Z}@>>> \tilde O(r,{\mathbb C})@>>> O(r,{\mathbb C})@>>> 1\\
@. @| @VVV @VVV @.\\
1@>>>{\mathbb Z}/2{\mathbb Z}@>>>
 \tilde {\bf O}(r)^{\rm an}@>>> {\bf O}(r)^{\rm an}@>>> 1
\end{CD}$$ 
of central extensions
of sheaves on $S^{\rm an}$.
\end{proof}

We prove the main result of
this section.

\begin{pr}\label{pr11}
Let $S$ be a normal scheme of finite type
over a noetherian ring over 
an algebraic closure
$\bar {\mathbb Q}$ of ${\mathbb Q}$.
Let $f:X\to S$ be a proper smooth morphism
of relative even dimension $n$.
Then, for every prime number $\ell$,
we have
$$sw_2(H^n_\ell(X/S))=
hw_2(H^n_{dR}(X/S))$$
in $H^2(S,{\mathbb Z}/2{\mathbb Z})$.
\end{pr}

If $K$ is an extension of
$\bar {\mathbb Q}$,
the remaining terms
in (\ref{eqcn})
are 0.
Hence Proposition \ref{pr11}
implies the assertion 4
in Theorem \ref{thm1}.

\begin{proof}
By a standard argument,
we may assume $S$ is of finite type over
$\bar {\mathbb Q}$.
By Lefschetz principle
$H^2(S,{\mathbb Z}/2{\mathbb Z})
\simeq
H^2(S_{\mathbb C},{\mathbb Z}/2{\mathbb Z})$,
it is reduced to the case where
$S$ is of finite type over
${\mathbb C}$.

We identify 
$H^2(S,{\mathbb Z}/2{\mathbb Z})
=H^2(S^{\rm an},{\mathbb Z}/2{\mathbb Z})$
by the canonical isomorphism.
Let $V$ be
the orthogonal local system
$R^nf^{\rm an}_*{\mathbb C}$
and
$D$ be
the symmetric bundle
$R^nf^{\rm an}_*\Omega^{\bullet {\rm an}}
_{X^{\rm an}/S^{\rm an}}$
on $X^{\rm an}$.
Then, we have
$sw_2(H^n_\ell(X/S))=
sw_2(V)$
and
$hw_2(H^n_{dR}(X/S))
=hw_2(D)$.
Since
$D=V\otimes_{\mathbb C}
{\mathcal  O}_{X^{\rm an}}$,
the assertion follows from
Lemma \ref{lm1a}.
\end{proof}

\section{Arithmetic argument}\label{sart}

We prove the assertion 5
of Theorem \ref{thm1}
by a global arithmetic
argument.
We will apply the following
to the moduli space
of hypersurfaces.

\begin{lm}\label{lmPS}
Let $S$ be a connected
normal noetherian scheme,
$K$ be the
fraction field of $S$
and $\bar K$
be an algebraic
closure of $K$.
Let $p:P\to S$ be a proper
smooth and
geometrically connected
scheme over $S$
and $U\subset P$
be an open subscheme.
Let $m$ be an integer
invertible on $S$.
We assume that
the following conditions
{\rm (1)--(3)}
are satisfied:

{\rm (1)} 
The open subscheme
$U\subset P$
is the complement of 
a divisor  $D\subset P$
flat over $S$.
Let $D^\circ\subset D$ be the largest
open subscheme 
smooth over $S$.
Then,  for every $s\in S$,
the fiber
$D^\circ_s$ is dense
in $D_s$.

{\rm (2)} 
For every irreducible
component $C$
of the geometric
generic fiber $D_{\bar K}$,
the divisor class
$[C]\in {\rm Pic}
(P_{\bar K})$
is in the image
of the multiplication
$m\times\colon 
{\rm Pic}
(P_{\bar K})
\to 
{\rm Pic}
(P_{\bar K})$.

{\rm (3)} 
$R^1p_*\mu_m=0$.

Let $D_1,\ldots,D_r$
be the irreducible components
of $D$ 
with the reduced closed subscheme
structures,
$F_i$ be the
fraction field of $D_i$
and $K_i$ be the local field
at the generic point of $D_i$
for $1\le i\le r$.
Then, the compositions
$H^2(U,\mu_m)
\to
H^2(K_i,\mu_m)
\to
H^1(F_i,\mu_m)$
with the boundary map
define an exact sequence
$$\begin{CD}
H^2(S,\mu_m)
@>>> 
H^2(U,\mu_m)
@>>>
H^2(U_{\bar K},\mu_m)
\oplus
\bigoplus_{i=1}^r
H^1(F_i,{\mathbb Z}/m{\mathbb Z}).
\end{CD}$$
\end{lm}

\begin{proof}
We consider the commutative diagram
$$\begin{CD}
H^2(S,\mu_m)@. @. \\
@VVV @. @.\\
H^2(P,\mu_m)
@>>>
H^2(U,\mu_m)
@>>>
\bigoplus_{i=1}^r
H^1(F_i,{\mathbb Z}/m{\mathbb Z})\\
@VVV @VVV @.\\
H^2(P_{\bar K},\mu_m)
@>>>
H^2(U_{\bar K},\mu_m).
\end{CD}$$
We show that
the conditions (1)--(3)
respectively
imply the exactness
of the middle row,
the injectivity
of the bottom horizontal
arrow
and the exactness of the
left column.
The assertion follows from this
by diagram chasing.

We show the exactness
of the middle row
assuming the condition (1).
The union $V=U\cup D^\circ$
is open in $P$ and we regard it 
as an open 
subscheme of $P$.
Since $D^\circ$ is a smooth divisor of $V$,
we have an exact sequence
$H^2(V,\mu_m)
\to 
H^2(U,\mu_m)
\to 
H^1(D^\circ,{\mathbb Z}/m{\mathbb Z})$
by relative purity.
Since 
$H^1(D^\circ,{\mathbb Z}/m{\mathbb Z})
\to 
\bigoplus_{i=1}^r
H^1(F_i,{\mathbb Z}/m{\mathbb Z})$
is injective,
it is sufficient to show that
the restriction map
$H^q(P,\mu_m)
\to 
H^q(V,\mu_m)$ is an isomorphism
for $q\le 2$.
Let $v:V\to S$ be
the restriction of $p:P\to S$.
By the assumption (1),
the codimension of the complement
$P_s\setminus V_s$ in each fiber $P_s$
is at least 2.
Hence the map $R^qv_!{\mathbb Z}/m{\mathbb Z}
\to R^qp_*{\mathbb Z}/m{\mathbb Z}$
is an isomorphism for 
$q\ge 2N-2$ where $N$ is the
relative dimension of $P$
over $S$.
By Poincar\'e duality,
the map $R^qp_*{\mathbb Z}/m{\mathbb Z}
\to R^qv_*{\mathbb Z}/m{\mathbb Z}$
is an isomorphism for 
$q\le 2$.
Hence the assertion follows
by comparing the Leray spectral sequences
$H^p(S,R^qp_*\mu_m)
\Rightarrow
H^{p+q}(P,\mu_m)$
and 
$H^p(S,R^qv_*\mu_m)
\Rightarrow
H^{p+q}(V,\mu_m)$.

Next we show that
the restriction map
$H^2(P_{\bar K},\mu_m)
\to 
H^2(U_{\bar K},\mu_m)$
is injective assuming the condition (2).
Let $C_1,\ldots,C_{r'}$
be the irreducible components
of $D_{\bar K}$.
By the same argument as above,
we obtain an exact sequence
$\bigoplus_{j=1}^{r'}
H^0(C_j,{\mathbb Z}/m{\mathbb Z})
\to 
H^2(P_{\bar K},\mu_m)
\to 
H^2(U_{\bar K},\mu_m).$
By the exact sequence
${\rm Pic}
(P_{\bar K})
\overset{m\times}
\to
{\rm Pic}
(P_{\bar K})
\to H^2(
P_{\bar K},\mu_m)$,
the injectivity is
in fact equivalent to
the condition (2).

Finally,
the condition (3)
implies the exactness
of the left column by the
Leray spectral sequence.
Thus the assertion is proved.
\end{proof}

Let $n\ge 0$ be an
integer and $d\ge 1$ be
an integer.
Let $P={\mathbb P}(
\Gamma({\mathbb P}^{n+1}_{\mathbb Z},
{\mathcal  O}(d))^\vee)$ be
the moduli space of
hypersurfaces of degree $d$
and $f:X\to P$ be
the universal family of
hypersurfaces.
More explicitly,
it is described as follows.
Let $T_0,\ldots,T_{n+1}
\in \Gamma({\mathbb P}^{n+1}_{\mathbb Z},
{\mathcal  O}(1))$
be the homogeneous coordinate
of ${\mathbb P}^{n+1}_{\mathbb Z}$.
The monomials
$T_0^{i_0}\cdots T_{n+1}^{i_{n+1}}$ 
of degree $i_0+\cdots +i_{n+1}=d$
form
a basis of 
$\Gamma({\mathbb P}^{n+1}_{\mathbb Z},
{\mathcal  O}(d))$.
Let $(A_{i_0,\ldots,i_{n+1}})$
for $i_0+\cdots +i_{n+1}=d$
be the dual basis of
$
\Gamma({\mathbb P}^{n+1}_{\mathbb Z},
{\mathcal  O}(d))^\vee$.
They form a basis of
$\Gamma(P,{\mathcal  O}(1))$
and $X$ is defined in 
${\mathbb P}^{n+1}_{\mathbb Z}\times
P$ by the equation
$F=0$
where $$F=\sum_{i_0+\cdots +i_{n+1}=d}
A_{i_0,\ldots,i_{n+1}}
T_0^{i_0}\cdots T_{n+1}^{i_{n+1}}$$
$\in \Gamma({\mathbb P}^{n+1}_{\mathbb Z}\times
P,{\mathcal  O}(d,1))$
is the universal polynomial.

Let $U\subset P$
be the open subscheme where
the universal hypersurface 
$f:X\to P$ is smooth.
If $X_K\subset {\mathbb P}^{n+1}_K$ 
is a smooth hypersurface
defined by a polynomial 

\noindent
$\sum_{i_0+\cdots +i_{n+1}=d}
a_{i_0,\ldots,i_{n+1}}
T_0^{i_0}\cdots T_{n+1}^{i_{n+1}}$
of degree
$d$ of coefficients in $K$,
it is the pull-back of
the universal family $X\to P$
by the map 
$\Spec K\to U$
defined by the coefficients
$(a_{i_0,\ldots,i_{n+1}})$.

\begin{lm}\label{lmDel}
Let $\Delta\subset X$ 
be the complement
of the largest
open subscheme of $X$
smooth over $P$.
We regard 
$\Delta$ as a closed
subscheme of $X$
defined by the equations
$\frac{\partial F}{\partial T_0}=
\cdots=\frac{\partial F}
{\partial T_{n+1}}=0$
and we put
$N+1={\rm rank}\ 
\Gamma({\mathbb P}^{n+1}_{\mathbb Z},
{\mathcal  O}(d))$.

{\rm 1.} 
The scheme $X$
is a ${\mathbb P}^{N-1}$-bundle
over ${\mathbb P}^{n+1}
_{\mathbb Z}$
with respect to the
first projection.
Consequently, it is regular
and irreducible.

{\rm 2.} 
The scheme $\Delta$
is a ${\mathbb P}^{N-(n+2)}$-bundle
over ${\mathbb P}^{n+1}
_{\mathbb Z}$ 
with respect to the
first projection.
Consequently,
it is regular
and irreducible.
The immersion
$\Delta\to X$
is a regular immersion
of codimension $n+1$.

{\rm 3.} 
Locally on $X$,
the coherent ${\mathcal  O}_X$-module
$\Omega^1_{X/P}$
is isomorphic to
the cokernel
of the map 
${\mathcal  O}_X
\to
{\mathcal  O}_X^{n+1}$
defined by a system
of generators
$a_1,\ldots,a_{n+1}$
of the ideal
${\mathcal  I}_\Delta
\subset {\mathcal  O}_X$.
\end{lm}

\begin{proof}
It suffices to
consider over 
the open subscheme
$D(T_i)\subset 
{\mathbb P}^{n+1}
_{\mathbb Z}$
for each $i=0,\ldots,n+1$.

1.
Since $X$ is defined
by the linear form
$F$ on $A_{i_0,\ldots,i_{n+1}}$,
the assertion follows.

2.
This is proved in 
\cite[Proposition 2.8]{dd}
using the fact that 
the closed scheme $\Delta
\subset 
P\times_{{\mathbb P}^{n+1}
_{\mathbb Z}}D(T_i)$
is defined by
$F=
\frac{\partial F}{\partial T_0}
=\cdots=
\frac{\partial F}{\partial T_{i-1}}
=
\frac{\partial F}{\partial T_{i+1}}
=\cdots=
\frac{\partial F}{\partial T_{n+1}}
=0$ that follows from
$d\cdot F=\sum_j 
T_j\cdot
\frac{\partial F}{\partial T_j}$.

3.
We define a regular
function $f_i$
on $X\times
_{{\mathbb P}^{n+1}
_{\mathbb Z}}D(T_i)$
by 
$F=f_i\cdot T_i^d$
and put $t_j=T_j/T_i$
Then, on $X\times
_{{\mathbb P}^{n+1}
_{\mathbb Z}}D(T_i)$,
the cokernel
of the map 
${\mathcal  O}_X
\to
{\mathcal  O}_X^{n+1}$
defined by 
$\frac{\partial f_i}
{\partial t_0},
\ldots,
\frac{\partial f_i}
{\partial t_{i-1}},
\frac{\partial f_i}
{\partial t_{i+1}},$
$\ldots,
\frac{\partial f_i}
{\partial t_{n+1}}$.
It follows from 
the proof of 2.\
that they form
a system of generators
of the ideal
${\mathcal  I}_\Delta
\subset {\mathcal  O}_X$.
\end{proof}

We regard
the image $D\subset P$ 
of $\Delta\subset X$
as a reduced closed subscheme.
On the scheme $D$,
the following is proved in
\cite[Lemma 2.10,Proposition 2.12)]{dd}.
The degree is computed in 
\cite[n$^{\rm o}$ 5]{De}.

\begin{lm}\label{lmDe}
The scheme $D$ 
is a divisor of degree
$(n+2)(d-1)^{n+1}$
flat over ${\rm Spec}\ {\mathbb Z}$
and has geometrically
irreducible fibers.
The largest
open subscheme of $D$
smooth over ${\rm Spec}\ {\mathbb Z}$
has non-empty intersection
with every fiber.
\end{lm}

\begin{pr}\label{prhp}
For a prime number $\ell>n+1$,
Conjecture {\rm \ref{cn2}}
is true for the universal family
$f_{U[\frac 1{2\ell}]}:
X_{U[\frac 1{2\ell}]}\to 
U[\frac 1{2\ell}]$
of smooth hypersurfaces.
\end{pr}

As Conjecture {\rm \ref{cn2}}
implies
Conjecture {\rm \ref{cn}},
Proposition
\ref{prhp}
implies
the assertion 5
of Theorem \ref{thm1}.
It also implies
the assertion 2
of Theorem \ref{thm1}
in the case $p=\ell=2$
since in this case we have
$n=0$.

\begin{proof}
Let $\delta
\in H^2(U[\frac 1{2\ell}],
{\mathbb Z}/2{\mathbb Z})$
be the difference
of the both sides of
(\ref{eqcncxS}).
If $\ell\neq 2$,
the difference $\delta$
lies in the image of
$H^2(U[\frac 12],
{\mathbb Z}/2{\mathbb Z})
\to H^2(U[\frac 1{2\ell}],
{\mathbb Z}/2{\mathbb Z})$,
by Lemma \ref{lmlS} and
Corollary \ref{cor3}.
For $\ell=2$,
it is obvious.
In order to show that
$\delta$
lies in the image of
$H^2({\mathbb Z}[\frac 12],
{\mathbb Z}/2{\mathbb Z})
\to H^2(U[\frac 1{2\ell}],
{\mathbb Z}/2{\mathbb Z})$,
we verify that
$U[\frac 12]
\subset P[\frac 12]\to 
{\rm Spec}\ {\mathbb Z}[\frac 12]$
and $m=2$ satisfy 
the conditions (1)--(3) of Lemma \ref{lmPS}.

Let $\xi$
denote the generic point
of the irreducible 
closed divisor $D\subset P$.
The conditions (1) and (2)
follow from Lemma \ref{lmDe}.
Since $P[\frac 12]$ is a projective space
over ${\rm Spec}\ {\mathbb Z}[\frac 12]$,
the condition (3) is 
also satisfied.
Thus, by Lemma \ref{lmPS},
the sequence
\setcounter{equation}0
\begin{equation}
H^2({\rm Spec}\ {\mathbb Z}[\frac 12],
{\mathbb Z}/2{\mathbb Z})
\to
H^2(U[\frac 12],{\mathbb Z}/2{\mathbb Z})
\to
H^2(U_{\bar {\mathbb Q}},{\mathbb Z}/2{\mathbb Z})
\oplus
H^1(\kappa(\xi),{\mathbb Z}/2{\mathbb Z})
\label{eqar}
\end{equation}
is exact.

Let $K$ be the completion
of the fraction field
of the discrete valuation
ring ${\mathcal  O}_{P,\xi}$.
Then, since $X$ is regular,
the base change
of $X\to P$
to the valuation ring 
${\mathcal  O}_K$
is also regular.
Further its closed fiber
has at most an ordinary
double point as
singularity
by \cite[Proposition 3.2]{SGA7}
(cf.\ Proof of \cite[Theorem 3.5]{dd}).
Hence 
by Corollary \ref{corodp},
the image
of $\delta$
in $H^1(\kappa(\xi),{\mathbb Z}/2{\mathbb Z})$
is 0.
By Proposition \ref{pr11},
its image in 
$H^2(U_{\bar {\mathbb Q}},
{\mathbb Z}/2{\mathbb Z})$
is also 0.
Thus, by the exact sequence
(\ref{eqar}),
it is in the image of
$H^2({\mathbb Z}[\frac12],
{\mathbb Z}/2{\mathbb Z})$.

Let $x\in U({\mathbb R})$
be the point corresponding
to the Fermat hypersurface.
By Proposition \ref{prHdg},
we also have
$x^*\delta=0$
in
$H^2({\mathbb R},
{\mathbb Z}/2{\mathbb Z})$.
Since the composite map
$$\begin{CD}
H^2({\mathbb Z}[\frac 12],
{\mathbb Z}/2{\mathbb Z})
@>>>
H^2(U[\frac1{2\ell}],{\mathbb Z}/2{\mathbb Z})
@>{x^*}>>
H^2({\mathbb R},
{\mathbb Z}/2{\mathbb Z})
\end{CD}$$
is an isomorphism,
we have
$\delta= 0$.
\end{proof}

\end{document}